\documentclass{article}
\title{Lindel\"of indestructibility, topological games and selection principles} 
\author{Marion Scheepers and Franklin D. Tall\footnote{Research supported by Grant A-7354 of the Natural Sciences and Engineering Research Council of Canada}}

\usepackage{amsfonts}
\newtheorem{theorem}{{\bf Theorem}}
\newtheorem{proposition}[theorem]{{\bf Proposition}}
\newtheorem{lemma}[theorem]{{\bf Lemma}}
\newtheorem{corollary}[theorem]{{\bf Corollary}}

\newcommand{\epf}{\Box\vspace{0.15in}}
\newtheorem{problem}{{\bf Problem}}

\newcommand{\naturals}{{\mathbb N}}
\newcommand{\integers}{{\mathbb Z}}
\newcommand{\reals}{{\mathbb R}}

\newcommand{\hilbertcube}{{\mathbb H}}

\newcommand{\forces}{\mathrel{\|}\joinrel\mathrel{-}}

\newcommand{\sone}{{\sf S}_1}
\newcommand{\gone}{{\sf G}_1}
\newcommand{\sfin}{{\sf S}_{fin}}
\newcommand{\gfin}{{\sf G}_{fin}}

\newcommand{\open}{\mathcal{O}}
\date{August 23, 2009}
\begin{document}
\maketitle
\begin{abstract} 
Arhangel'skii proved that if a first countable Hausdorff space is Lindel\"of, then its cardinality is at most $2^{\aleph_0}$. Such a clean upper bound for Lindel\"of spaces in the larger class of spaces whose points are ${\sf G}_{\delta}$ has been more elusive. In this paper we continue the agenda started in \cite{FDT}, of considering the cardinality problem for spaces satisfying stronger versions of the Lindel\"of property. Infinite games and selection principles, especially the Rothberger property, are essential tools in our investigations\footnote{\lowercase{{\bf {\uppercase{K}}ey words and phrases:} {\uppercase{I}}ndestructibly {\uppercase{L}}indel\"of, infinite game, {\uppercase{R}}othberger space, {\uppercase{G}}erlits-{\uppercase{N}}agy space, {\uppercase{M}}enger property, {\uppercase{H}}urewicz property, measurable cardinal, weakly compact cardinal, real-valued measurable cardinal, 
{\uppercase{C}}ohen reals, random reals.}}. 
\end{abstract}

A topological space is Lindel\"of if each open cover contains a countable subset that covers the space. 
Alexandrov asked if in the class of first countable Hausdorff spaces, every Lindel\"of space has cardinality at most $2^{\aleph_0}$. Arhangel'skii \cite{arh} proved that the answer is ``yes". This focuses attention on the larger class of spaces in which ``points are ${\sf G}_{\delta}$" - \emph{i.e.}, the class of topological spaces in which each point is an intersection of countably many open sets. Such spaces are ${\sf T}_1$ but not necessarily Hausdorff. Arhangel'skii showed that in the class of spaces with points ${\sf G}_{\delta}$ each Lindel\"of space has cardinality less than the least measurable cardinal. Juh\'asz \cite{Juhasz} showed that this bound is sharp: There are such Lindel\"of spaces of arbitrary large cardinality below the least measurable cardinal. Juh\'asz's examples are not Hausdorff spaces and the cardinality of the underlying spaces has countable cofinality. Shelah \cite{SS} showed that no Lindel\"of space with points ${\sf G}_{\delta}$ can be of weakly compact cardinality. Gorelic \cite{IG} showed that it is relatively consistent that the Continuum Hypothesis (CH) holds, that $2^{\aleph_1}$ is arbitrarily large, and there is a zero-dimensional regular Lindel\"of space with points ${\sf G}_{\delta}$ and of cardinality $2^{\aleph_1}$. This improved earlier results of Shelah \cite{SS}, Shelah-Stanley \cite{SSt} and Velleman \cite{Ve} which showed that either by countably closed forcing, or by assuming V=L, one could obtain such a space of cardinality $\aleph_2$, consistent with CH. Little else is known about cardinality and the Lindel\"of property in the class of spaces with points ${\sf G}_{\delta}$\footnote{R. Knight has published a paper \cite{Knight} claiming that V = L implies there are Lindelof spaces with points ${\sf G}_\delta$ of size $\aleph_n$, for each $n \le \omega$ (See Theorems 3.6.1, 3.6.2 and 3.6.3 of \cite{Knight}). Questions have been raised as to whether the proof is correct. In \cite{FDT} a similar result was attributed to C. Morgan, but that proof turned out to be incorrect.}. It is not even known (absolutely or consistently) if all such ${\sf T}_2$ spaces have cardinality $\le 2^{\aleph_1}$. 

This cardinality problem is of interest also for stronger versions of the Lindel\"of property. An interesting strengthening of the Lindel\"of property is identified and studied in \cite{FDT}: Call a Lindel\"of space \emph{indestructibly Lindel\"of} if in any generic extension by a countably closed partially ordered set the space is still Lindel\"of. \cite{FDT} Theorem 19 shows that Juh\'asz's examples are destructible. And \cite{FDT} Theorem 4 shows, assuming large cardinals, that it is consistent that in the class of spaces with points ${\sf G}_{\delta}$, each indestructibly Lindel\"of space has cardinality at most $2^{\aleph_0}$. 
It is also shown in \cite{FDT} (in the paragraph following Theorem 18) that Gorelic's space is indestructibly Lindel\"of. Thus, the statement that indestructibly Lindel\"of spaces with points ${\sf G}_{\delta}$ are of cardinality $\le 2^{\aleph_0}$ is (modulo large cardinals) independent. This independence can also be obtained in another way if the separation hypothesis is weakened: A. Dow \cite{Adow} showed that adding $\aleph_1$ Cohen reals converts every ground model Lindel\"of space to an indestructibly Lindel\"of space in the generic extension. This gives the consistency of the existence of non-${\sf T}_2$ indestructibly Lindel\"of spaces with points ${\sf G}_{\delta}$ of arbitrary large cardinality below the first measurable cardinal: Juh\'asz's spaces from the ground model retain their cardinality but acquire Lindel\"of indestructibility, and measurable cardinals from the ground model remain measurable.

\cite{FDT} Theorem 3 gives a combinatorial characterization of indestructibly Lindel\"of. In Section 1, starting from this combinatorial characterization, we show how to characterize indestructibly Lindel\"of game-theoretically (Theorem \ref{FDT}) and examine the determinacy of this game. This analysis leads us to two natural strengthenings of indestructibly Lindel\"of. For one of these strengthenings, all spaces with points ${\sf G}_{\delta}$ and this strengthening have cardinality $\le 2^{\aleph_0}$ (Theorem \ref{GT}). The classical selection property introduced by Rothberger is the other natural strengthening (Corollary \ref{rothbconnect}). These topics are illustrated with some examples collected in Section 4.

In Section 2 we focus on the Rothberger property. 
In \cite{DTW2}, using the technique of n-dowments, it was shown that a variety of non-covering and non-generalized-metric properties were preserved by Cohen reals. Since these arguments mimic measure-theoretic ones, the same preservation arguments work for adding random reals. This work was extended by \cite{GJT} and by \cite{Iwasa}. It turns out that covering properties are also often preserved. From this previous work, it was not at all evident that there would be some covering or non-covering properties which would be preserved by one of these kinds of extensions, but not the other. Here we shall show that while ``non-Rothberger" is not preserved by Cohen extensions (of uncountable cardinality) it is preserved by random real extensions. 
Specifically: We show that forcing with $\kappa$ Cohen reals ($\kappa$ uncountable) converts all ground model Lindel\"of spaces to spaces with the Rothberger property (Theorem \ref{cohenrothberger}). This improves the above-mentioned result of A. Dow. We show that forcing with $\kappa$ random reals preserves the Rothberger property, and for a large class of spaces, forcing with $\kappa$ random  reals preserves the non-Rothberger property (Theorems \ref{randomrothbergerpreserve} and \ref{randomrothberger}). We also show that Rothberger spaces not only retain the Lindel\"of property under countably closed forcing (since they are indestructibly Lindel\"of) but in fact they retain the Rothberger property (Theorem \ref{countablyclosed}). 

These three results lead to several illuminating pieces of information about the (insufficiently studied) Rothberger property in general spaces: Every Rothberger space with points ${\sf G}_{\delta}$ has cardinality no larger than the least real-valued measurable cardinal (Theorem \ref{rvmrothberger}). If it is consistent that there is a measurable cardinal, then it is consistent that all Rothberger spaces with points ${\sf G}_{\delta}$ have cardinality at most $2^{\aleph_0}$. 
(Even ${\sf T}_3$) Rothberger spaces of cardinality larger than $\aleph_1$ need not contain a Lindel\"of subspace of cardinality $\aleph_1$. However it is consistent, assuming the consistency of the existence of a measurable cardinal, that Rothberger ${\sf T}_2$ spaces of character less than or equal to $\aleph_1$ include Rothberger subspaces of size at most $\aleph_1$ (Theorem \ref{collapserothb}). We also find a new forcing proof for Shelah's result that there is no points ${\sf G}_{\delta}$ Lindel\"of topology on a set of weakly compact cardinality (Proposition \ref{noweakcptlindelof}).

In Section 3 we examine a covering property introduced by Gerlits and Nagy \cite{GN}. The Gerlits-Nagy property is closely related to the Rothberger property. Indeed, for regular spaces it is undecidable by ZFC whether these two properties are different. Study of the Gerlits-Nagy property requires considering the preservation, by forcing, of selection properties introduced by Hurewicz. One of these selection properties is called the Menger property, and the other the Hurewicz property. Unlike the Menger property and the Rothberger property, both Cohen real and random real forcing preserve being not Hurewicz (Theorems \ref{nonhurewiczpreserve} and \ref{randomnonhurewicz}). While Cohen real forcing does not preserve the Hurewicz property in general, random real forcing does (Theorem \ref{nonhurewicz}). We also found a class of separable metrizable spaces (which include Sierpi\'nski sets) for which Cohen forcing does preserve the Hurewicz property (Theorem \ref{stronghurewiczpreserve}). Thus it is consistent that the square of a set of real numbers with the Gerlits-Nagy property need not have the Gerlits-Nagy property (Corollary \ref{problem6.6}). This solves Problem 6.6 of \cite{Ts}.

In Section 4 we remark on a few examples in light of our results. We show that Gorelic's example has both the Rothberger and Hurewicz properties while Juh\'asz's examples have the Hurewicz property but not the Rothberger property.

{\flushleft{\bf Acknowledgement:}} We thank the referee for a very careful reading of the paper, and for pointing out that Cohen real forcing does not in general preserve the Hurewicz property.

\section{Indestructibly Lindel\"of and an infinite game.}

Let $\mathcal{A}$ and $\mathcal{B}$ be collections of sets and let $\alpha$ be an ordinal number. Then the game $\gone^{\alpha}(\mathcal{A},\mathcal{B})$ between players ONE and TWO is played as follows: They play an inning for each $\beta<\alpha$. In the $\beta$-th inning ONE first chooses some $\mathcal{U}_{\beta}\in\mathcal{A}$; TWO responds by choosing some set $T_{\beta}\in \mathcal{U}_{\beta}$. A play
\[
  \mathcal{U}_0,\, T_0,\, \cdots,\, \mathcal{U}_{\beta},\, T_{\beta},\, \cdots : \beta<\alpha
\]
is won by TWO if $\{T_{\beta}:\beta<\alpha\}$ is a member of $\mathcal{B}$; else, ONE wins.

Let $\open$ denote the set of all open covers of space $X$. In \cite{FG} Galvin introduced the version $\gone^{\omega}(\open,\open)$ of the above game. 

A collection $(U_f: f\in \bigcup_{\alpha<\omega_1} \,^{\alpha}\omega)$ of open subsets of a space $X$ is said to be a \emph{covering tree} if for each $\alpha < \omega_1$, for each $f\in\,^\alpha\omega$, the set $\{U_{f\frown\{(\alpha,n)\}}:n<\omega\}$ is a cover of $X$. In Theorem 3 of \cite{FDT} it is proved that a space $X$ is indestructibly Lindel\"of if, and only if, it is Lindel\"of, and for each covering tree $(U_f: f\in \bigcup_{\alpha<\omega_1} \,^{\alpha}\omega)$ of $X$ there is an $f\in \bigcup_{\alpha<\omega_1}\, ^{\alpha}\omega$, such that the set $\{ U_{f\lceil_{\beta}}: \beta < domain(f)\}$ is a cover of $X$.

\begin{theorem}\label{FDT}
For a Lindel\"of space $X$ the following are equivalent:
\begin{enumerate}
  \item{X is indestructibly Lindel\"of.}
  \item{ONE has no winning strategy in ${\sf G}^{\omega_1}_1(\mathcal{O},\mathcal{O})$.}
\end{enumerate}
\end{theorem}
{\bf Proof:} To see that $(1)\Rightarrow (2)$, observe that if $F$ is a strategy for player ONE then there is a natural covering tree associated with $F$: We may assume, since $X$ is Lindel\"of, that in each inning $F$ calls on ONE to play a countable open cover of $X$. Thus enumerate ONE's first move $F(\emptyset)$ as $(U_{(0,n)}:n<\omega)$. When TWO chooses $U_{(0,n_0)}\in F(\emptyset)$, ONE's response $F(U_{(0,n_0)})$ may be enumerated as $(U_{\{(0,n_0),\, (1,n)\}}:n<\omega)$. When TWO chooses $U_{\{(0,n_0),\, (1,n_1)\}} \in F(U_{\{(0,n_0)\}})$, ONE's response $F(U_{(0,n_0)},\, U_{\{(0,n_0),\, (1,n_1))\}})$ is enumerated as $(U_{\{(0,n_0),\, (1,n_1), (2,n))\}}:n<\omega)$, and so on.

This produces the covering tree $(U_f: f\in \bigcup_{\alpha<\omega_1} \,^{\alpha}\omega)$ of $X$, associated as above with ONE's strategy $F$. Now apply the fact that $X$ is indestructibly Lindel\"of to find an $\alpha<\omega_1$ and an $f\in \, ^{\alpha}\omega$ such that  $\{ U_{f\lceil_{\beta}}: \beta < \alpha\}$ is a cover of $X$. TWO wins the play
\[
  F(\emptyset),\, U_{f\lceil_1},\, F(U_{f\lceil_1}), \cdots,\, F(U_{f\lceil_{\beta}}: \beta < \gamma), \, U_{f\lceil_{\gamma+1}},\, \cdots,
\]
demonstrating that $F$ is not a winning strategy for ONE.

To see that $(2)\Rightarrow (1)$, observe that a covering tree rather directly defines a strategy of ONE in the game $\gone^{\omega_1}(\open,\open)$, and that a winning play of TWO against this strategy witnesses the indestructibility of $X$ for this covering tree. $\epf$

Theorem \ref{FDT} suggests a number of at least formal strengthenings of the notion of indestructibly Lindel\"of. We consider specifically the following two: One strengthening is to require that TWO has a winning strategy in $\gone^{\omega_1}(\open,\open)$; the second is to require that ONE has no winning strategy already in the shorter game $\gone^{\omega}(\open,\open)$.

\begin{center}{\bf When TWO has a winning strategy}\end{center}

If TWO has a winning strategy in $\gone^{\omega_1}(\open,\open)$ then every open cover of the space contains a subcover of cardinality at most $\aleph_1$. In Theorem \ref{GT} we do not assume that the space is Lindel\"of and use an idea of F. Galvin \cite{FG} in the proof. 
\begin{theorem}\label{GT} Let $X$ be a space in which each point is a ${\sf G}_{\delta}$. If TWO has a winning strategy in the game $\gone^{\omega_1}(\open,\open)$ then $\vert X\vert \le 2^{\aleph_0}$.
\end{theorem}
{\bf Proof:} 
Since each point of $X$ is a {\sf G}$_{\delta}$, fix for each $x\in X$ a sequence $(U_n(x):n<\omega)$ of open sets with $\{x\} = \bigcap_{n<\omega}U_n(x)$. Let $F$ be a strategy for TWO. Let $\alpha<\omega_1$ as well as a sequence $(\mathcal{U}_{\beta}:\beta<\alpha)$ of open covers of $X$ be given.

{\flushleft{\bf Claim 1:}} There is an $x\in X$ such that for each open set $U\subseteq X$ with $x\in U$, there is an open cover $\mathcal{U}$ of $X$ with 
\[
  U = F((\mathcal{U}_{\beta}:\beta<\alpha)\frown(\mathcal{U})).
\]
{\bf Proof of Claim 1:} If not, choose for each $x\in X$ an open set $U_x\subseteq X$ with $x\in U_x$ and $U_x$ is a counterexample to the above statement. Put $\mathcal{U} = \{U_x:x\in X\}$. Then $\mathcal{U}$ is an open cover of $X$, and $F((\mathcal{U}_{\beta}:\beta<\alpha)\frown(\mathcal{U}))\in\mathcal{U}$, contradicitng the definition of members of $\mathcal{U}$.

By the claim choose $x_{\emptyset}\in X$ so that there is for each open set $U$ with $x_{\emptyset}\in U$ an open cover $\mathcal{U}$ of $X$ with $U = F(\mathcal{U})$. For each $n$ choose an open cover $\mathcal{U}_{\{(0,n)\}}$ of $X$ with $U_n(x_{\emptyset})=F(\mathcal{U}_{\{(0,n)\}})$.

For each $n<\omega$ choose $x_{\langle n\rangle}$ for open cover $\mathcal{U}_{\{(0,n)\}}$ as in the claim. Then, for each $n_0$ and for each $n$ choose an open cover $\mathcal{U}_{\{(0,n_0),(1,n)\}}$ of $X$ such that 
\[
  U_n(x_{\langle n_0\rangle}) = F(\mathcal{U}_{\{(0,n_0)\}},\mathcal{U}_{\{(0,n_0),(1,n)\}}).
\] 
Then for each $n_0,\, n_1<\omega$ choose $x_{\langle n_0,\, n_1\rangle}\in X$ as in the claim for the sequence  $(\mathcal{U}_{\{(0,n_0)\}},\mathcal{U}_{\{(0,n_0),(1,n)\}})$ of open covers. For each $n<\omega$ choose an open cover $\mathcal{U}_{\{(0,n_0),\, (1,n_1),\, (2,n)\}}$ of $X$ such that
\[
  U_n(x_{\langle n_0,\, n_1\rangle}) = F(\mathcal{U}_{\{(0,n_0)\}},\mathcal{U}_{\{(0,n_0),(1,n)\}},\, \mathcal{U}_{\{(0,n_0),\, (1,n_1),\, (2,n)\}}),
\] 
and so on. In general, let $\alpha<\omega_1$ be given and assume we have selected for each $f\in\bigcup_{\beta<\alpha}\, ^{\beta}\omega$ points $x_f\in X$ and open covers $\mathcal{U}_f$ of $X$ such that for all $\beta<\alpha$ and for $f\in\, ^{\beta}\omega$ and $n< \omega$,
\[
  U_n(x_f) = F((\mathcal{U}_{f\lceil_{\gamma}}:\gamma<\beta)\frown(\mathcal{U}_{f\cup\{(\beta,n)\}})).
\]
Consider $f\in\, ^{\alpha}\omega$. Applying Claim 1 to the sequence $(\mathcal{U}_{f\lceil_{\beta}}:\beta<\alpha)$ of open covers, choose an $x_f\in X$ as in Claim 1. Then for each $n<\omega$ choose an open cover $\mathcal{U}_{f\cup\{(\alpha,n)\}}$ of $X$ with
\[
  U_n(x_f) = F(\mathcal{U}_{f\lceil_{\beta}}:\beta<\alpha)\frown(\mathcal{U}_{f\cup\{(\alpha,n)\}})).
\]

This specifies how to recursively choose for each $f\in\bigcup_{\alpha<\omega_1}\, ^{\alpha}\omega$ a point $x_f\in X$ and an open cover $\mathcal{U}_f$ such that for each such $f$ and each $n<\omega$, if $f\in\, ^{\gamma}\omega$, then
\[
  U_n(x_f) = F(\mathcal{U}_{f\lceil_{\beta}}:\beta<\gamma)\frown(\mathcal{U}_{f\cup\{(\gamma,n)\}})).
\]
Note that $\{x_{f}:f\in\bigcup_{\alpha<\omega_1}\, ^{\alpha}\omega\}$ has cardinality at most $2^{\aleph_0}$. 

{\flushleft{\bf Claim 2:}} $X = \{x_{f}:f\in\bigcup_{\alpha<\omega_1}\, ^{\alpha}\omega\}$.\\
{\bf Proof of Claim 2:} For if not, choose $x\in X\setminus\{x_{f}:f\in\bigcup_{\alpha<\omega_1}\, ^{\alpha}\omega\}$. Choose $n_0$ with $x\not\in U_{n_0}(x_{\emptyset})$. ONE plays $\mathcal{U}_{\{(0,n_0)\}}$ and TWO responds using $F$. Choose $n_1$ with $x\not\in U_{n_1}(x_{\langle n_0\rangle})$. Next ONE plays $\mathcal{U}_{\{(0,n_0),(1,n_1)\}}$ and TWO responds using $F$, and so on. In this way we build an $F$-play which is lost by TWO since TWO never covers the point $x$. This contradicts the fact that $F$ is a winning strategy for TWO. $\epf$

A result of \cite{DG} can be formulated as follows.
\begin{theorem}[Daniels-Gruenhage]\label{DGth} If $X$ is hereditarily Lindel\"of then TWO has a winning strategy in $\gone^{\omega_1}(\open,\open)$. 
\end{theorem}

Combining Theorems \ref{GT} and \ref{DGth} gives a new proof for an old result of de Groot:
\begin{corollary} If $X$ is hereditarily Lindel\"of and each point is a ${\sf G}_{\delta}$-set, 
then $\vert X\vert\le 2^{\aleph_0}$. Thus if $X$ is hereditarily Lindel\"of and ${\sf T}_2$, then $\vert X\vert\le 2^{\aleph_0}$.
\end{corollary}

\begin{center}{\bf Determinacy of the game $\gone^{\omega_1}(\open,\open)$}\end{center}

If a topological space $X$ has cardinality $\le \aleph_1$, then TWO has a winning strategy in $\gone^{\omega_1}(\open,\open)$ on that space. Thus:

\begin{corollary}\label{chandtwo} Assume the Continuum Hypothesis. For a Lindel\"of space $X$ with points ${\sf G}_{\delta}$ the following are equivalent:
\begin{enumerate}
  \item{TWO has a winning strategy in $\gone^{\omega_1}(\open,\open)$.}
  \item{$\vert X\vert\le \aleph_1$.}
\end{enumerate} 
\end{corollary}

An infinite two-player game in which neither player has a winning strategy is said to be \emph{undetermined}. 
In the class of Lindel\"of spaces with points ${\sf G}_{\delta}$ the determinacy of the game $\gone^{\omega_1}(\open,\open)$ is (modulo large cardinals) independent of ZFC. This can be seen as follows:

\begin{theorem}[Gorelic, \cite{IG}]\label{gorelic}
For any cardinal number $\aleph_{\alpha}$ it is consistent, relative to the consistency of {\sf ZFC}, that $2^{\aleph_0} = \aleph_1$ and there is a zero-dimensional Lindel\"of space $X$ with points ${\sf G}_{\delta}$ and $\vert X\vert = 2^{\aleph_1}\ge \aleph_{\alpha}$.
\end{theorem}

Shelah's earlier proof of consistency of existence of such a space in a model where $\aleph_1 = 2^{\aleph_0} < 2^{\aleph_1} = \aleph_2$ finally appeared in \cite{SS}. Note that such a space is necessarily ${\sf T}_3$: points ${\sf G}_{\delta}$ implies ${\sf T}_1$, and ${\sf T}_1$ plus zero-dimensional implies ${\sf T}_3$. Tall showed in \cite{FDT} that both Shelah's and Gorelic's examples are indestructibly Lindel\"of. Theorems \ref{FDT} and \ref{GT} imply that the game $\gone^{\omega_1}(\open,\open)$ is undetermined on these two spaces.

In Theorem 4 of \cite{FDT} it is proved consistent, relative to the consistency of the existence of a supercompact cardinal, that the indestructibly Lindel\"of spaces with points ${\sf G}_{\delta}$ are of small cardinality:
\begin{theorem}[Tall]\label{consistency} 
If it is consistent that there is a supercompact cardinal, then it is consistent that the Continuum Hypothesis holds and every indestructibly Lindel\"of space with points ${\sf G}_{\delta}$ is of cardinality $\le 2^{\aleph_0}$.
\end{theorem}

This yields the following determinacy result:
\begin{corollary}\label{determinacy} If it is consistent that there is a supercompact cardinal, then it is consistent that the game $\gone^{\omega_1}(\open,\open)$ is determined on all Lindel\"of spaces with points ${\sf G}_{\delta}$.
\end{corollary}
{\bf Proof:} Consider a model of Theorem \ref{consistency}. CH holds in this model. Consider in this model a Lindel\"of space $X$ with points ${\sf G}_{\delta}$. If $\vert X\vert\le\aleph_1$ then by Corollary \ref{chandtwo} TWO has a winning strategy in $\gone^{\omega_1}(\open,\open)$. If $\vert X\vert>\aleph_1$, then by Theorem 4 of \cite{FDT}, $X$ is destructible. By Theorem \ref{FDT} ONE has a winning strategy in $\gone^{\omega_1}(\open,\open)$. $\epf$ 

\begin{center}{\bf When ONE does not win the shorter game $\gone^{\omega}(\open,\open)$}\end{center}

There is a natural selection principle associated with the game $\gone^{\alpha}(\mathcal{A},\mathcal{B})$: The symbol $\sone^{\alpha}(\mathcal{A},\mathcal{B})$ denotes the statement:
\begin{quote} For each sequence $(A_{\gamma}:\gamma<\alpha)$ of elements of $\mathcal{A}$, there is a sequence $(B_{\gamma}:\gamma<\alpha)$ such that for each $\gamma$ we have $B_{\gamma}\in A_{\gamma}$, and $\{B_{\gamma}:\gamma<\alpha\}\in\mathcal{B}$. 
\end{quote}
If $\kappa$ is an initial ordinal with the same cardinality as $\alpha$, then $\sone^{\kappa}(\mathcal{A},\mathcal{B})$ holds if, and only if, $\sone^{\alpha}(\mathcal{A},\mathcal{B})$ holds. Thus, we consider $\sone^{\alpha}(\mathcal{A},\mathcal{B})$ only for initial ordinals. 

If ONE does not have a winning strategy in $\gone^{\alpha}(\mathcal{A},\mathcal{B})$, then $\sone^{\alpha}(\mathcal{A},\mathcal{B})$ holds. Thus, if $X$ is indestructibly Lindel\"of, then it has the property $\sone^{\omega_1}(\open,\open)$. This observation gives an easy proof that the compact space $^{\omega_1}2$ is not indestructibly Lindel\"of: For $\alpha<\omega_1$ put $U^{\alpha}_i = \{f\in\,^{\omega_1}2: f(\alpha)=i\}$. Then each $\mathcal{U}_{\alpha}:=\{U^{\alpha}_0,\, U^{\alpha}_1\}$ is an open cover of $^{\omega_1}2$ and the sequence $(\mathcal{U}_{\alpha}:\alpha<\omega_1)$ witnesses that $\sone^{\omega_1}(\open,\open)$ fails for this space.
\begin{problem} When does $\sone^{\omega_1}(\open,\open)$ imply that ONE has no winning strategy in $\gone^{\omega_1}(\open,\open)$?
\end{problem}

Rothberger introduced the selection property $\sone^{\omega}(\open,\open)$ in \cite{rothberger38}. Spaces with the property $\sone^{\omega}(\open,\open)$ are said to have the Rothberger property, and are called \emph{Rothberger spaces}. We shall use the notation $\sone^{\omega}(\open,\open)$ and the terminology ``Rothberger space". Note that the Rothberger property is inherited by closed subspaces and preserved by countable unions and continuous images. We shall often use the following observation about spaces $X$ with the property $\sone^{\omega}(\open,\open)$: for each sequence $(\mathcal{U}_n:n<\omega)$ of open covers of $X$ there is a sequence $(U_n:n<\omega)$ such that for each $n$, $U_n\in\mathcal{U}_n$, and for each  $x\in X$ the set $\{n:x\in U_n\}$ is infinite. This follows by partitioning $\omega$ into countably many disjoint infinite sets $(S_n, \, n<\omega)$, and then applying $\sone^{\omega}(\open,\open)$ to each of the sequences $(\mathcal{U}_n:n\in S_m)$, $m<\omega$, of open covers of $X$.

\begin{theorem}[Pawlikowski, \cite{JP}]\label{pawlikowski}
A topological space $X$ has property $\sone^{\omega}(\open,\open)$ if, and only if, ONE has no winning strategy in the game $\gone^{\omega}(\open,\open)$.
\end{theorem}

If ONE has no winning strategy in $\gone^{\omega}(\open,\open)$, then ONE has no winning strategy in $\gone^{\omega_1}(\open,\open)$. It follows that:
\begin{corollary}\label{rothbconnect} Rothberger spaces are indestructibly Lindel\"of.
\end{corollary}

Unlike the case with the game $\gone^{\omega_1}(\open,\open)$ (see Corollary \ref{determinacy}), the game $\gone^{\omega}(\open,\open)$ is provably undetermined in the class of points ${\sf G}_{\delta}$ Lindel\"of spaces: J.T. Moore's L-space (Section 4, Example 2) and Todor\v{c}evic's stationary Aronszajn lines (Section 4, Example 6). These are uncountable Rothberger spaces with points ${\sf G}_{\delta}$. Since these spaces are Rothberger spaces, Pawlikowski's Theorem implies that ONE has no winning strategy in $\gone^{\omega}(\open,\open)$. In \cite{FG} Galvin proved that if TWO has a winning strategy in $\gone^{\omega}(\open,\open)$ in a points ${\sf G}_{\delta}$ Lindel\"of space $X$, then $X$ is countable. Since Moore's and Todor\v{c}evic's spaces are uncountable, also TWO has no winning strategy in $\gone^{\omega}(\open,\open)$. 

Examples 2 and 3 are spaces of cardinality $\aleph_1$, and the spaces of Example 6 can also be taken to be of cardinality $\aleph_1$. Thus, for each of these spaces TWO has a winning strategy in $\gone^{\omega_1}(\open,\open)$. 
In Section 4, Example 3 we examine Gorelic's space more closely and show that it is a Rothberger space, whence ONE has no winning strategy in $\gone^{\omega}(\open,\open)$ and TWO has no winning strategy in $\gone^{\omega_1}(\open,\open)$ on this example.

As will be seen later, in the class of metrizable spaces the determinacy of the game $\gone^{\omega}(\open,\open)$ is independent of ZFC.

\section{The Rothberger property}

The Rothberger property arose in connection with the classical Borel Conjecture. E. Borel \cite{Borel} defined a subset $X$ of the real line $\reals$ to be \emph{strong measure zero} if: For each sequence $(\epsilon_n:n\in\naturals)$ of positive real numbers there is a sequence $(I_n:n\in\naturals)$ of open intervals such that for each $n$ $length(I_n)\le \epsilon_n$, and $X\subseteq\bigcup_{n\in\naturals}I_n$. Borel conjectured that all strong measure zero sets of reals are countable. The statement that strong measure zero sets of real numbers are countable sets is known as the \emph{Borel Conjecture}. Rothberger \cite{rothberger38} showed that for metrizable spaces the Rothberger property implies strong measure zero. Indeed, for metrizable spaces the Rothberger property is equivalent to having strong measure zero in all equivalent metrics \cite{FM}. 

Rothberger spaces have been rarely considered outside the metric context. Indeed, to our knowledge, cardinality restrictions on Rothberger spaces have received no prior attention outside the metric context. 

\begin{center}{\bf Rothberger spaces and Cohen reals}\end{center}

For an infinite cardinal number $\kappa$, let ${\sf Fn}(\kappa,2)$ denote the set of finite binary functions with domain a finite subset of $\kappa$. For $f$ and $g$ elements of ${\sf Fn}(\kappa,2)$ write $f<g$ to denote that $f$ extends $g$. Thus, $({\sf Fn}(\kappa,2),<)$ is one of several partially ordered sets for adding $\kappa$ Cohen reals. ${\sf Fn}(\kappa,\omega)$, the set of functions with domain a finite subset of $\kappa$ and values in $\omega$, and with the same order as above, leads to the same generic extension. We shall use the symbol ${\mathbb P}(\kappa)$ to denote either of these posets for adding $\kappa$ Cohen reals.

\begin{theorem}\label{cohenrothberger}
Let $\kappa$ be an uncountable cardinal number. 
If $X$ is a Lindel\"of space, then ${\mathbf 1}_{{\mathbb P}(\kappa)} \forces ``{\check{X}} \mbox{ has property }\sone^{\omega}(\open,\open)"$.
\end{theorem}
{\bf Proof:} Let $X$ be a Lindel\"of space. Let $\kappa$ be an uncountable cardinal. Everywhere below when we write ${\mathbb P}(\kappa)$, we may consider this forcing as expressed as ${\mathbb P}(\kappa+\omega_1)$. It is known (see \cite{Adow} or \cite{GJT}) that ${{\mathbf 1}}_{{\mathbb P}(\kappa)} \forces ``{\check{X}} \mbox{ is Lindel\"of}"$. Consider a  ${\mathbb P}(\kappa)$-name for a sequence of open covers of $X$, say $(\dot{\mathcal{U}}_n:n\in\naturals)$. Since the ground-model topology of $X$ is a basis for the topology in the generic extension, and since $X$ is Lindel\"of, we may assume that 
\[
  {{\mathbf 1}}_{{\mathbb P}(\kappa)} \forces ``\mbox{Each $\dot{\mathcal{U}}_n$ is a countable cover by ground-model open sets}"
\] 
Indeed, for each $n$ and $m$ we can choose a maximal antichain $A^n_m$ of ${\mathbb P}(\kappa)$ such that the $m$-th element of $\dot{\mathcal{U}}_n$ is decided by this antichain. Considering ${\mathbb P}(\kappa)$ as ${\mathbb P}(\kappa+\omega_1)$, and since each member of each $A^n_m$ has domain a finite subset of $\kappa$, there is an $\alpha<\kappa$ such that in fact 
\[
  {{\mathbf 1}}_{{\mathbb P}(\alpha)} \forces ``\mbox{Each $\dot{\mathcal{U}}_n$ is a countable cover by ground-model open sets}".
\] 
Let $G$ be ${\mathbb P}(\kappa)$-generic. Then $G_{\alpha} = \{p\lceil_{\alpha}: p\in G\}$ is ${\mathbb P}(\alpha)$-generic over V. In V[G$_{\alpha}$] write $\mathcal{U}_n=(U^n_m:m\in\naturals)$. Define, for each $x\in X$ the function $f_x\in \,^{\naturals}\naturals$ by $f_x(n) = \min\{m:\, x\in U^n_m\}$. Now $\{f_x:x\in X\}\in$ V[G$_{\alpha}$]. If c $\in$ V[G] is any Cohen real over V[G$_{\alpha}$], then for each x in X, $\{n\in\naturals: c(n) = f_x(n)\}$ is infinite, and this $(U^n_{c(n)}:n\in\naturals)$ witnesses $\sone^{\omega}(\open,\open)$ for $X$ in V[G].
$\epf$

Thus, forcing with a sufficient number of Cohen reals preserves being a Rothberger space, but fails to preserve not being Rothberger. 

\begin{problem}\label{singlecohentorothberger} Is it consistent that adjoining a single Cohen real converts each ground-model Lindel\"of space to a Rothberger space or at least makes it indestructible? 
\end{problem}

\begin{corollary}[A. Dow]\label{dowth} Adjoining $\aleph_1$ Cohen reals renders each ground-model Lindel\"of space indestructible.
\end{corollary}
{\bf Proof:} By Theorem \ref{cohenrothberger} each ground-model Lindel\"of space acquires the Rothberger property. By Theorem \ref{pawlikowski} in the extension ONE has no winning strategy in the game $\gone^{\omega_1}(\open,\open)$. Apply Theorem \ref{FDT}. $\epf$\\
Dow's proof in \cite{Adow} uses a quite non-trivial elementary submodel argument. 

\begin{center}{\bf Rothberger spaces and random reals}\end{center}

Let ${\mathbb B}(\kappa)$ be the partially ordered set for adding $\kappa$ random reals. As was noted in the remarks in and following the proof of Corollary 2.5 of \cite{GJT}, forcing with ${\mathbb B}(\kappa)$ preserves the Lindel\"of property. We now show: Forcing with ${\mathbb B}(\kappa)$ preserves being Rothberger. Contrary to the case of Cohen reals, for a large class of Lindel\"of spaces random reals preserve not being Rothberger. 

\begin{lemma}\label{randomcover} Let $X$ be a Lindel\"of space and let $\dot{\mathcal{U}}$ be a ${\mathbb B}(\kappa)$ name such that 
\[
  {\bf 1}_{{\mathbb B}(\kappa)}\forces``\dot{\mathcal{U}} \mbox{ is an open cover of }\check{X}."
\]
Then for each $x\in X$ and each $k<\infty$ there is a neighborhood $N_k(x)$ of $x$ such that
\[
  \mu(\Vert(\exists U\in\dot{\mathcal{U}})(\check{N}_k(x)\subseteq U) \Vert)>1-\frac{1}{2^{k+1}}
\]
\end{lemma}
{\bf Proof:} Since $X$ is Lindel\"of and random reals preserve Lindel\"of, and since the ground-model open sets form a basis of the topology $X$ has in the generic extension, we may assume that
\[
  {\bf 1}_{{\mathbb B}(\kappa)}\forces``\dot{\mathcal{U}} \mbox{ is a countable open cover of }\check{X} \mbox{ by ground-model sets}"
\]
and choose ${\mathbb B}(\kappa)$ names $\dot{U}_n$, $n<\infty$ such that 
\[
  {\bf 1}_{{\mathbb B}(\kappa)}\forces``\dot{\mathcal{U}} =\{\dot{U}_n:\,n<\omega\} \mbox{ and each }\dot{U}_n \mbox{ is in the ground-model}". 
\]
We may assume each $\dot{U}_n$ is of the form $\{(\check{U}^n_k,b^n_k):k<\omega\}$ where for each $n$ and each $k$ $U^n_k$ is  a ground model open set  and $\{b^n_k:k<\infty\}$ is a maximal antichain of ${\mathbb B}(\kappa)$.

Now consider $x\in X$. For each $n$ and $k$ define 
\[
  V^n_k(x) = \left\{\begin{tabular}{ll}
                        $X$ & $\mbox{if }x\not\in U^n_k$\\
                        $U^n_k$ & $\mbox{ otherwise}$
                     \end{tabular}
            \right.   
\]
For each $k$ define $M_k(x) = \cap_{i\le k,\\ j\le k} V^i_j(x)$. Note that 
\[
  c_k = \Vert(\exists i\le k)( \check{M}_k(x)\subseteq \dot{U}_i)\Vert \ge \sup\{b^i_j:\, i,\, j\le k \mbox{ and }x\in U^i_j\}.
\]
For $k<\ell$ we have $\mu(c_k)\le \mu(c_{\ell})$, and since $\mu(\Vert\dot{\mathcal{U}}\mbox{ is a cover}\Vert)=1$ we have $\lim_{k\rightarrow\infty}\mu(c_k)=1$. For each $m$ choose $k_m$ so large that $\mu(c_{k_m})>1-\frac{1}{2^{m+1}}$ and put $N_m(x) = M_{k_m}(x)$.
$\epf$

\begin{theorem}\label{randomrothbergerpreserve} If $X$ is a Rothberger space then ${\bf 1}_{{\mathbb B}(\kappa)} \forces ``\check{X} \mbox{ is a Rothberger space}"$.
\end{theorem}
{\bf Proof:} Choose ${\mathbb B}(\kappa)$ names $\dot{\mathcal{U}}_n$ such that 
\[
  {\bf 1}_{{\mathbb B}(\kappa)}\forces ``(\dot{\mathcal{U}}_n:n<\infty) \mbox{ is a sequence of open covers of }\check{X}".
\]

By Lemma \ref{randomcover} for each $n$, and for each $x\in X$ choose a neighborhood $W_n(x)$ such that 
$\mu(\Vert(\exists U\in\dot{\mathcal{U}}_n)(\check{W}_n(x)\subseteq U) \Vert)>1-\frac{1}{2^{n+1}}$. Then $\mathcal{W}_n=\{W_n(x):x\in X\}$ is an open cover of $X$. In the ground model apply the fact that $X$ is Rothberger to the sequence $(\mathcal{W}_n:n<\infty)$ of open covers of $X$. For each $n$ we choose a set $S_n\in\mathcal{W}_n$ such that for each $x\in X$ there are infinitely many $n$ with $x\in S_n$.

Observe that we have for each $n$: $\mu(\Vert(\exists U\in\dot{\mathcal{U}}_n)(\check{S}_n\subseteq U) \Vert)>1-\frac{1}{2^{n+1}}$. Thus, choose by the Fullness Lemma (see \cite{jech},  Lemma 18.6) for each $n$ a ${\mathbb B}(\kappa)$ name $\dot{U}_n$ such that $\mu(\Vert \dot{U}_n\in\dot{\mathcal{U}}_n \mbox{ and } \check{S}_n\subseteq \dot{U}_n\Vert) \ge 1-\frac{1}{2^{n+1}}$. 

{\flushleft{{\bf Claim:}}} ${\bf 1}_{{\mathbb B}(\kappa)}\forces`` \check{X}\subseteq\bigcup_{n<\infty}\dot{U}_n"$.

For suppose not and choose $b$ so that $b\forces``\check{X}\not\subseteq\bigcup_{n<\infty}\dot{U}_n"$. Choose $c<b$ and $x\in X$ so that $c\forces ``\check{x}\not\in\bigcup_{n<\infty}\dot{U}_n"$. Then choose $n$ so large that $\mu(c)>\frac{1}{2^n}$. Choose a $k$ larger than $n$ with $x\in S_k$. Now define
\[
  d = c\wedge\Vert ``\dot{U}_k\in\dot{\mathcal{U}}_k \mbox{ and } \check{S}_k\subseteq \dot{U}_k"\Vert.
\]
We have $\mu(d)>0$ and since $d\le c$ we have $d\forces ``\check{x}\not\in\bigcup_{n<\infty}\dot{U}_n"$. But also $x\in S_k$ and $d\forces ``\check{S}_k\subseteq \dot{U}_k"$, a contradiction. The proof is complete.
$\epf$

Next we examine preservation of non-Rothberger: Let $n$ be a positive integer. A family $\mathcal{A}$ of subsets of a set $X$ has \emph{degree} $n$ if $n$ is the smallest positive integer such that for each $x\in X$ the set $\{A\in\mathcal{A}:x\in A\}$ has cardinality at most $n$. We shall say that a topological space is $n$-\emph{dimensional} if $n$ is the smallest nonnegative integer such that each open cover of the space has an open refinement of degree $n+1$. And a space is said to be \emph{strongly countable dimensional} if it is a union of countably many closed subsets, each of finite dimension. This notion was introduced by Nagata \cite{nagata} and Smirnov \cite{smirnov}, independently.

If a strongly countable dimensional space does not have the Rothberger property, then some closed subset of finite dimension does not have the Rothberger property. This follows from the fact that the Rothberger property is inherited by closed subspaces, and since a union of countably many subspaces, each with the Rothberger property, again has the Rothberger property. We now show that for strongly countable dimensional ground model spaces random reals preserve the property of not being a Rothberger space. In the process we use the following fact - see ``equation" (\ref{positive1}) in the proof of Theorem \ref{randomrothberger} below: If a space $X$ has the Rothberger property and if a sequence $(\mathcal{U}_n:n<\infty)$ of open covers is given, then there is a sequence $(U_n:n<\infty)$ such that for each $n$ we have $U_n\in\mathcal{U}_n$ and for each $x\in X$ there are infinitely many $n$ with $x\in U_n$:  First write the natural numbers as a union of countably many disjoint infinite sets $Y_k,\, k<\infty$, and apply the Rothberger property to each subsequence $(\mathcal{O}_n:n\in Y_k)$.

\begin{lemma}\label{orderlemma} If $\dot{U}$ is a ${\mathbb B}(\kappa)$-name and A and B are ground model sets
such that 
\[
  \mu(\Vert\dot{U} = {\check{A}}\Vert \wedge \Vert\dot{U} = {\check{B}}\Vert)>0
\]
then $A = B$.
\end{lemma}
{\bf Proof:} 
Consider a generic filter $G$ which contains $\Vert\dot{U} = {\check{A}} \Vert \wedge \Vert\dot{U} = {\check{B}}\Vert$. In the generic extension, $A = {\check{A}}_G = \dot{U}_G = {\check{B}}_G = B$. But $A$ and $B$ are ground model sets. $\epf$

\begin{theorem}\label{randomrothberger}
Let $\kappa$ be a cardinal number with uncountable cofinality. 
If $X$ is a strongly countable dimensional Lindel\"of space and if
\[
  {{\mathbf 1}}_{{\mathbb B}(\kappa)} \forces ``{\check{X}} \mbox{ has the Rothberger property}",
\]
 then $X$ has the Rothberger property in the ground-model.
\end{theorem}
${\bf Proof}$ Suppose that ${{\mathbf 1}}_{{\mathbb B}(\kappa)} \forces ``{\check{X}} \mbox{ has the Rothberger property}"$ while, in fact, $X$ does not have the Rothberger property. Choose a positive integer $m$ and a closed subset $Y$ of $X$ which is $(m-1)$-dimensional, but does not have the Rothberger property. Let $(\open_n:n<\infty)$ be a sequence of open covers of $Y$ witnessing this. We may assume each $\open_n$ has degree $m$, and that each $\open_{n+1}$ refines $\open_n$. Enumerate each $\open_n$ bijectively as $(U^n_k:k\in I_n)$, where $I_n$ is either an initial segment of the set of natural numbers, or the set of natural numbers.  
Note that
\[
  {{\mathbf 1}}_{{\mathbb B}(\kappa)} \forces ``{\check{Y}} \mbox{ is closed in ${\check{X}}$, thus has the Rothberger property }."
\]

For each $n$, put $\epsilon_n = \frac{1}{(n+1)^2\dot m}$. Define $\psi:\naturals\rightarrow\naturals$ as follows:
$\psi(1) = \frac{1}{\epsilon_1}$ and for each $n$, $\psi(n+1) = \psi(n) + \frac{1}{\epsilon_{n+1}}$.

In $V^{\mathbb B(\kappa)}$ look at $({\check{\open}}_{\psi(n)}: n<\infty)$, a sequence of open covers of $Y$. Since ${{\mathbf 1}}_{{\mathbb B}(\kappa)} \forces ``{\check{Y}} \mbox{ has the Rothberger property}"$, choose a sequence $(\dot{V}_n:n<\infty)$ of ${\mathbb B} (\kappa)$-names such that 
\begin{equation}\label{positive1}
  {{\mathbf 1}}_{{\mathbb B}(\kappa)} \forces ``(\forall n)(\emptyset\neq \dot{V}_n \in {\check{\open}}_{\psi(n)}) \mbox{ and }(\forall x\in {\check{Y}})(\exists^{\infty}_n)(x\in \dot{V}_n)".
\end{equation}

For each $n$ define
\[
  \mathcal{W}_n:=\{O^{\psi(n)}_j\in\open_{\psi(n)}: \mu(\| \dot{V}_n = {\check{O}}^{\psi(n)}_j\|) >\epsilon_n\}.
\]
Since elements of $\open_{\psi(n)}$ are listed without repetition, Lemma \ref{orderlemma} implies that $C_n = \vert\mathcal{W}_n\vert \le m\cdot(n+1)^2$. List the elements of each $\mathcal{W}_n$ as $(O^{\psi(n)}_{m^n_j},\, 1\le j \le C_n)$.

Using the fact that for each $n$, $\mathcal{O}_{n+1}$ refines $\mathcal{O}_n$, we can choose $O^n_{m_n}\in\open_n,\, n< \infty$, such that there is for each $U\in\bigcup_{n<\infty}\mathcal{W}_n$ an $n$ with $U\subseteq O^n_{m_n}$. 

Since $(\open_n:n<\infty)$ witnessed that $Y$ does not have the Rothberger property, choose 
\[
  x \in Y\setminus \bigcup_{n<\infty}O^n_{m_n}.
\]

{\flushleft{\bf Claim 1:}} For each $n$, $\mu(\| {\check{x}}\in \dot{V}_n\|) \le \frac{1}{(n+1)^2}$.\\
{\bf Proof of the Claim:} If not, choose $n$ with $\mu(\| {\check{x}}\in \dot{V}_n\|) > \frac{1}{(n+1)^2}$. Since $\open_{\psi(n)}$ covers $Y$, the set $I = \{j: x\in  O^{\psi(n)}_j\}$ is nonempty and for each $j\in I$ we have $\mu(\|\check{x}\in \check{O}^{\psi(n)}_j\|)=1$, and so $\mu(\|\check{x}\in \dot{V}_n \cap (\cap_{j\in I} \check{O}^{\psi(n)}_j)\|)> \frac{1}{(n+1)^2}$. But then as $\mu(\|\dot{V}_n\in{\check{\open}_{\psi(n)}}\|) = 1$ we have 
\[
  \frac{1}{(n+1)^2} < \mu(\Vert \dot{V}_n \in \{\check{O}^{\psi(n)}_j:j\in I\}\Vert)
\]
 and $\open_{\psi(n)}$ has degree $m$, for some $j\in I$ we have $\mu(\|\dot{V}_n = \check{O}^{\psi(n)}_j \|)>\frac{1}{m\cdot (n+1)^2}$ and so for this $j$, $O^{\psi(n)}_j \in\mathcal{W}_n$, contradicting the choice of $x$. Claim 1 is proven.

Consider any $b\in {\mathbb B}(\kappa)$ with $\mu(b)>0$, and choose $n$ so large that $\mu(b) > \sum_{j=n}^{\infty}\frac{1}{j^2}$.
 With Claim 1 proven, define
\[
  a = b \setminus(\bigvee_{k\ge n}\|\check{x}\in \dot{V}_k\|).
\] 
We have $a< b$, $\mu(a)>0$ and 
\begin{equation}\label{negative1}
  a \forces ``(\forall k\ge n)(\check{x}\not\in \dot{V}_k)".
\end{equation}
But then (\ref{negative1}) contradicts (\ref{positive1}). Thus, $X$ has the Rothberger property also. $\epf$

In our first proof of Theorem \ref{randomrothberger} we used the hypothesis that the space is zero-dimensional. Subsequently we improved the hypothesis that the space be zero-dimensional to the current hypothesis of strong countable dimensionality. We were not able to eliminate the dimension-theoretic hypothesis altogether from Theorem \ref{randomrothberger}:
\begin{problem} Is the dimension hypothesis in Theorem \ref{randomrothberger} necessary?
\end{problem}

A large class of topological spaces is covered by the zero-dimensional case of Theorem \ref{randomrothberger}:
\begin{lemma}\label{TychonoffRothb}
${\sf T}_{3}$ spaces with the property $\sone^{\omega}(\open,\open)$ are zero-dimensional.
\end{lemma}
{\bf Proof:} If $X$ is ${\sf T}_3$ and Lindel\"of, it is ${\sf T}_4$ and thus ${\sf T}_{3\frac{1}{2}}$. Since $X$ is {\sf T}$_{3\frac{1}{2}}$ fix an embedding of $X$ into a power of the unit interval $I$, say $F:X\longrightarrow I^{\kappa}$. For each $\alpha\in \kappa$ the projection $\pi_{\alpha}$ onto the $\alpha$-th coordinate is continuous, and so $\pi_{\alpha}\circ F[X]\subseteq[0,1]$ has property $\sone^{\omega}(\open,\open)$, and thus is zero-dimensional. But then $\Pi_{\alpha<\kappa} \pi_{\alpha}\circ F[X]$ is zero-dimensional. Since the latter set contains $F[X]$, it follows that $F[X]$, and thus $X$, is zero-dimensional.
$\epf$\\

It is easy to convert, by forcing, certain ground model spaces to spaces which have dimension zero in the generic extension. In the generic extension by one Cohen real every ground model separable metric space acquires dimension zero. To see this, first observe that the real line of the ground model acquires Lebesgue measure zero: Temporarily let $\mu$ denote Lebesgue measure and let $\mu^*$ denote the outer measure. In the ground model choose for each $\epsilon>0$  a sequence $(\epsilon_n:n<\omega)$ of positive reals such that $\sum_{n<\omega}\epsilon_n <\epsilon$. We show that 
\[
  {\mathbf 1}_{{\mathbb P}(\omega)}\forces ``\mu(\check{\reals}) = 0".
\]  
If not, choose $p\in{\mathbb P}(\omega)$ and a positive integer $n$ such that $p\forces``\mu^*(\check{\reals})>\frac{1}{\check{n}}$".
Consider in the ground model the sequence $(\epsilon_m:m<\omega)$ associated with $\frac{1}{n}=\epsilon$, and for each $m$ define $\mathcal{U}_m =\{J^m_k: k<\omega\}$ where $\{J^m_k:k<\omega\}$ enumerates the set of intervals of length less than $\epsilon_m$ which have rational endpoints. For each $x\in \reals$  define $f_x$ so that for all $m$, $f_x(m)=\min\{k:x\in J^m_k\}$. Note that each $f_x$ is in the ground model. For each $x\in \reals$ and $m>\vert p\vert$ the set $D^x_m =\{q\in {\mathbb P}(\omega): (\exists k>m)(f_x(k)=q(k)\}$ is dense in ${\mathbb P}(\omega)$, and thus dense below $p$. Then any generic filter containing $p$ meets each $D^x_m$, and so the Cohen real obtained from this generic filter, say its ${\mathbb P}(\omega)$-name is $\dot{f}$, is infinitely often equal to each $f_x$. This implies that $p\forces ``\check{\reals}\subseteq \cup_{n<\omega}\dot{J}^n_{\dot{f}(n)}$". But since $\sum_{m<\omega}\epsilon_m < \frac{1}{n}$, it follows that $p\forces ``\mu^*(\check{\reals})<\frac{1}{\check{n}}$", a contradiction. 

Since Lebesgue measure zero sets of real numbers are zero-dimensional, it follows that the ground model set of reals is zero dimensional in the generic extension. But then the Hilbert cube of the ground model is in the extension a product of zero dimensional spaces and so zero dimensional. Separable metric spaces from the ground model embed homeomorphically into the Hilbert cube of the ground model, thus in the generic extension are subspaces of a zero dimensional space and so themselves are zero dimensional.

This specific consequence of Cohen real forcing holds more generally. 
\begin{corollary}\label{lindeloftozerodim}
If $X$ is a ${\sf T}_{3}$ Lindel\"of space and $\kappa$ is any cardinal number, then ${{\mathbf 1}}_{{\mathbb P}(\kappa)} \forces ``\check{X} \mbox{ is zero-dimensional}"$.
\end{corollary}

Similar ideas give the following result for spaces not necessarily Lindel\"of:
\begin{corollary}\label{arbitrarytozerodim}
Let $\kappa$ be a cardinal number.  In some ground model $M$, let $X$ is a ${\sf T}_{3\frac{1}{2}}$ space with $\vert X\vert<\kappa$. Consider a generic extension $N$ of $M$ such that $N\models \vert\reals\vert\ge \kappa$. Then $N \models ``{X} \mbox{ is zero-dimensional}"$.
\end{corollary}
{\bf Proof:} As in Lemma \ref{TychonoffRothb} we see that in $N$ each projection $\Pi_{\gamma}\circ F[X]\subseteq \reals$ has cardinality less than $2^{\aleph_0}$ and thus is zero-dimensional. $\epf$

\begin{center}{\bf Rothberger spaces and countably closed forcing.}\end{center}

Since Rothberger spaces are indestructibly Lindel\"of, they remain Lindel\"of under countably closed forcing. The following result shows a little more:

\begin{theorem}\label{countablyclosed} Let $({\mathbb P},<)$ be a countably closed partially ordered set. For $X$ a space the following are equivalent:
\begin{enumerate}
  \item {$X$ is a Rothberger space.}
  \item {${\bf 1}_{{\mathbb P}}\forces ``\check{X} \mbox{ is a Rothberger space}$."}
\end{enumerate}
\end{theorem}
{\bf Proof:} $1\Rightarrow 2$: Let $(X,\mathcal{T})$ be a space which has property $\sone^{\omega}(\open,\open)$. Let $(\dot{\mathcal{U}}_n:n\in\naturals)$ be a name for a sequence of open covers of $X$. That is,
\[
  1_{{\mathbb P}}\forces ``(\dot{\mathcal{U}}_n:n\in\naturals) \mbox{ is a sequence of open covers for } \check{X}."
\]
Let $G$ be ${\mathbb P}$-generic and in $V[G]$ consider the sequence of open covers. Since $X$ was, in the ground-model, $\sone^{\omega}(\open,\open)$, it is in $V[G]$ still Lindel\"of. We may also assume that each $\mathcal{U}_n$ consists of sets from the ground-model topology $\mathcal{T}$ on $X$. Thus, for each $n$ there is in $V[G]$ a function $f_n:\omega\longrightarrow\mathcal{T}$ such that $\mathcal{U}_n = f_n[\omega]$. 
Since the forcing is countably closed, each $f_n$ is a member of the ground model, and thus each $\mathcal{U}_n$ is a member of the ground model. But then apply $\sone^{\omega}(\open,\open)$ in the ground model to this sequence to obtain $U_n\in\mathcal{U}_n$, $n<\infty$, such that $\{U_n:n<\infty\}$ is a cover of $X$. 

{\flushleft{$2\Rightarrow 1$:}} Assume that ${\bf 1}_{{\mathbb P}}\forces ``\check{X} \mbox{ is a Rothberger space}$". Suppose that contrary to 1, in the ground model $X$ is not a Rothberger space. Choose in the ground model a sequence $(\mathcal{U}_n:n<\infty)$ of open covers of $X$ witnessing that $X$ is not Rothberger. By 2, choose a ${\mathbb P}$-name $\dot{f}$ such that 
\[
   {\bf 1}_{{\mathbb P}}\forces ``\dot{f}:\check{\omega}\longrightarrow \bigcup_{n<\infty}\check{\mathcal{U}}_n\mbox{ and }
  (\forall n)(\dot{f}(n)\in\check{\mathcal{U}}_n) \mbox{ and }\{\dot{f}(n):n\in\check{\omega}\} \mbox{ covers }\check{X}".
\] 
Choose for each $n$ a $p_n\in{\mathbb P}$, and a $U_n\in \mathcal{U}_n$ such that $p_{n+1}<p_n$, and $p_n\forces``\dot{f}(\check{n}) = \check{U}_n"$. Since ${\mathbb P}$ is countably closed choose $p\in{\mathbb P}$ such that for all $n$ we have $p<p_n$. Then $p\forces``(\forall n)(\check{U}_n\in\check{\mathcal{U}}_n \mbox{ and } \check{X} = \bigcup_{n<\infty}\check{U}_n)"$. But all parameters in the statement forced by $p$ are in the ground model; thus the statement is true in the ground model. This contradicts that the sequence $(\mathcal{U}_n:n<\infty)$ of open covers of $X$ witnesses that $X$ is not Rothberger.
$\epf$

\begin{center}{\bf Cardinality upper bounds on points ${\sf G}_{\delta}$ Rothberger spaces}\end{center}

It is consistent that all metrizable Rothberger spaces are countable. This follows from the Borel Conjecture. 
Carlson \cite{CarlsonBC} showed that the Borel Conjecture is equivalent to the statement that all strong measure zero metric spaces are countable. Thus the Borel Conjecture implies that all metrizable Rothberger spaces are countable. Sierpi\'nski \cite{sierpinski} proved that the Continuum Hypothesis implies the negation of the Borel Conjecture. R. Laver \cite{laver} proved that the Borel Conjecture is consistent relative to the consistency of ZFC only. In Laver's model, $2^{\aleph_0} = \aleph_2$. 

\begin{proposition}[Judah, Shelah, Woodin]\label{realspoor}
For each cardinal number $\kappa$ it is consistent that $2^{\aleph_0}>\kappa$ and every metrizable Rothberger space is countable.
\end{proposition}
{\bf Proof:} Consider a model of Borel's Conjecture. By Theorem \ref{randomrothberger} we have that upon forcing with ${\mathbb B}(\kappa)$ every ground-model set of real numbers with the Rothberger property is countable. It is not known if uncountable metric Rothberger spaces can be introduced by random real forcing. But by a result of Judah, Shelah and Woodin \cite{JSW} if the ground model is Laver's model, then subsequent forcing with ${\mathbb B}(\kappa)$ introduces no uncountable metrizable Rothberger spaces. 
$\epf$

\begin{corollary}\label{rvmandsize} If it is consistent that there is a real-valued measurable cardinal then it is consistent that each metrizable Rothberger space is countable and there is a real-valued measurable cardinal $\le 2^{\aleph_0}$.
\end{corollary}
{\bf Proof:} Assume it is consistent that there is a measurable cardinal, $\kappa$. Then it is consistent that $\kappa$ is measurable and CH holds (collapse the continuum to $\aleph_1$ with countable conditions, if necessary, and apply \cite{LS}, Theorem 3: ``Mild" forcing preserves measurability.). Then use Laver's method to force Borel's Conjecture. By \cite{LS}, Theorem 3, $\kappa$ still is measurable. Now force with ${\mathbb B}(\kappa)$. As noted in the proof of Theorem \ref{realspoor}, in the generic extension all metrizable Rothberger spaces are countable. By a theorem of Solovay $\kappa=2^{\aleph_0}$ is real-valued measurable.
$\epf$

Arhangel'skii's result trivially implies that if $\kappa$ is the least measurable cardinal, then each points ${\sf G}_{\delta}$ Rothberger spaces is of cardinality $<\kappa$. We shall see in Proposition \ref{measopt} that this upper bound cannot be lowered in ZFC. But a better theorem can be proved for the points ${\sf G}_{\delta}$ Rothberger spaces:
\begin{theorem}\label{rvmrothberger} If a space has the Rothberger property and each point is a ${\sf G}_{\delta}$ then its cardinality is less than the smallest real-valued measurable cardinal.
\end{theorem}
{\bf Proof:} Suppose that $\kappa$ is a real-valued measurable cardinal and that $X$ is a Lindel\"of space with all points ${\sf G}_{\delta}$ and $\vert X\vert\ge \kappa$. Choose a subset $Y$ of $X$ with $\vert Y\vert=\kappa$, and let $\mu:\mathcal{P}(Y)\rightarrow \lbrack 0, 1\rbrack$ be a countably additive atomless measure with $\mu(Y)=1$. For each $x\in X$ choose a sequence $(U_n(x):n<\infty)$ of neighborhoods with $\{x\} = \bigcap_{n<\infty}U_n(x)$ and $U_{n+1}(x)\subseteq U_n(x)$, all $n$.

For each $m$, choose for each $x$ an $n(m,x)$ such that $\mu(Y\cap U_{n(m,x)}(x))< \frac{1}{2^{m+2}}$. Then
\[
  \mathcal{U}_m = \{U_{n(m,x)}(x):x\in X\}
\]
is an open cover of $X$. For the sequence $(\mathcal{U}_m:m<\infty)$ of open covers of $X$, consider any sequence $(U_m:m<\infty)$ with $U_m\in\mathcal{U}_m$ for each $m$. Then we have for each $m$ that $\mu(Y\cap U_m)<\frac{1}{2^{m+2}}$, and so
\[
  \mu(Y\cap(\bigcup_{m<\infty}U_m))\le \sum_{m=1}^{\infty}\frac{1}{2^{m+2}} \le \frac{1}{2} < \mu(Y).
\]
But then the sequence $(U_m:m<\infty)$ does not cover $Y$ and so not $X$. It follows that $X$ is not Rothberger.
$\epf$\\ 
The bound given in Theorem \ref{rvmrothberger} cannot be lowered in ZFC: 

\begin{theorem}\label{rvmoptimality} If it is consistent that there is a real-valued measurable cardinal, then it is consistent that there is a real-valued measurable cardinal $\kappa\le 2^{\aleph_0}$ and for each cardinal $\lambda<\kappa$ there is a space $X$ with points ${\sf G}_{\delta}$ and the Rothberger property, with $\lambda\le \vert X\vert<\kappa$. 
\end{theorem}
{\bf Proof:} If it is consistent that there is a real-valued measurable cardinal, then it is consistent that there is a measurable cardinal \cite{Srvm}. Thus, assume that $\kappa$ is the least measurable cardinal. There are Lindel\"of spaces with points ${\sf G}_{\delta}$ of arbitrary large cardinality below $\kappa$, for example Juh\'asz's spaces descibed in Section 4, Example 4. First force with $\aleph_1$ Cohen reals. Then by Theorem \ref{cohenrothberger} the ground-model versions of Juh\'asz's spaces are Rothberger spaces with points ${\sf G}_{\delta}$ of arbitrary large cardinality below $\kappa$. By Theorem 3 of \cite{LS}, $\kappa$ is still measurable in this generic extension. Now force with ${\mathbb B}(\kappa)$ and apply Solovay's theorem \cite{Srvm} that ${\bf 1}_{{\mathbb B}(\kappa)}\forces ``2^{\aleph_0} =\check{\kappa} \mbox{ is real-valued measurable}"$.

In the resulting model $2^{\aleph_0}$ is real-valued measurable and there are Rothberger spaces with points ${\sf G}_{\delta}$ of arbitrary large cardinality below $2^{\aleph_0}$.
$\epf$\\
Note that in the proof of Theorem \ref{rvmoptimality} we may also force with ${\mathbb B}(\nu)$ for $\nu>\kappa$, thus causing the least real-valued measurable cardinal to be less than $2^{\aleph_0}$. 

Theorem \ref{rvmrothberger} also provides an alternative proof to Arhangel'skii's theorem:
\begin{corollary}[Arhangel'skii] Every points ${\sf G}_{\delta}$ Lindel\"of space has cardinality less that the least (two-valued) measurable cardinal.
\end{corollary}
{\bf Proof:} Assume that on the contrary $\kappa$ is the least two-valued measurable cardinal and that $X$ is a points ${\sf G}_{\delta}$ space of cardinality at least $\kappa$. First force with $\aleph_1$ Cohen reals, converting $X$ to a Rothberger space while preserving the measurability of $\kappa$. Now force with $\kappa$ random reals, converting $\kappa$ into a real-valued measurable cardinal while preserving the Rothberger property of $X$. This gives a contradiction to the result of Theorem \ref{rvmrothberger}. $\epf$

Let us also note that in the absence of real-valued measurable cardinals, the least measurable cardinal may be the least upper bound on the cardinality of points ${\sf G}_{\delta}$ Rothberger spaces:

\begin{proposition}\label{measopt}
If it is consistent that there is a measurable cardinal, then it is consistent that there is a measurable cardinal $\kappa$ and for each $\lambda<\kappa$ there is a points ${\sf G}_{\delta}$ Rothberger space of cardinality at least $\lambda$.
\end{proposition}
{\bf Proof:} Assume that $\kappa$ is the least measurable cardinal. There are Lindel\"of spaces with points ${\sf G}_{\delta}$ of arbitrary large cardinality below $\kappa$, for example Juh\'asz's spaces descibed in Section 4, Example 4. First force with $\aleph_1$ Cohen reals. Then by Theorem \ref{cohenrothberger} the ground-model versions of Juh\'asz's spaces are Rothberger spaces with points ${\sf G}_{\delta}$ of arbitrary large cardinality below $\kappa$, and by Theorem 3 of \cite{LS}, $\kappa$ is still measurable in this generic extension. 
$\epf$

It is a trivial consequence of Tall's result, Theorem \ref{consistency}, that if it is consistent that there is a supercompact cardinal then it is consistent that all points ${\sf G}_{\delta}$ Rothberger spaces have cardinality $\le \aleph_1$ and $2^{\aleph_0}=\aleph_1$. However, we don't know if either of the following is possible (even assuming the consistency of appropriate large cardinals):
\begin{problem}\label{carddistinguishrothberger}
Is it possible that all points ${\sf G}_{\delta}$ Rothberger spaces have cardinality $\le\aleph_1$ while there is an indestructible Lindel\"of space with points ${\sf G}_{\delta}$ and of cardinality larger than $\aleph_1$?
\end{problem}

Proposition \ref{twomeasures} records that if it is consistent that there is a measurable cardinal, then it is consistent that every points ${\sf G}_{\delta}$ Rothberger space has cardinality less than $2^{\aleph_0}$ and there are points ${\sf G}_{\delta}$ Lindel\"of spaces of arbitrarily large cardinality below the first (two-valued) measurable cardinal.

\begin{proposition}\label{twomeasures} Assume that it is consistent that there is a measurable cardinal. Then the conjunction of the following list of statements is consistent:
\begin{enumerate}
  \item{Every Rothberger space with points ${\sf G}_{\delta}$ is of cardinality less than $2^{\aleph_0}$;}
  \item{For each cardinal $\lambda<2^{\aleph_0}$ there is a Rothberger space with points ${\sf G}_{\delta}$ such that $\lambda<\vert X\vert<2^{\aleph_0}$;}
  \item{For each cardinal $\lambda>2^{\aleph_0}$ which is less than the first measurable cardinal there is a points ${\sf G}_{\delta}$ Lindel\"of space $X$ such that $\lambda<\vert X\vert$.}
  \end{enumerate} 
\end{proposition} 

\begin{problem}\label{carddistinguishrothberger2}
Is it possible that all points ${\sf G}_{\delta}$ Rothberger spaces have cardinality $\le\aleph_1$ while $2^{\aleph_0}>\aleph_1$?
\end{problem}

Since the real line is indestructibly Lindel\"of, a positive answer to Problem \ref{carddistinguishrothberger2} gives a positive answer to Problem \ref{carddistinguishrothberger}.

\begin{center}{\bf When Lindel\"of implies Rothberger.}\end{center} 

In Theorem \ref{PspaceHurewicz} we shall see a topological hypothesis under which Lindel\"of implies Rothberger (and more). Here we explore set theoretic hypotheses that imply the equivalence of the Lindel\"of and Rothberger properties: Let $\mathcal{M}_{\reals}$ denote the collection of first category subsets of the real line. Then ${\sf cov}(\mathcal{M}_{\reals})$ denotes the least cardinality of a family of first category sets whose union covers the real line.
\begin{proposition}\label{covmeager} For an infinite cardinal number $\kappa$ the following are equivalent:
\begin{enumerate}
  \item{Every Lindel\"of space of cardinality at most $\kappa$ has the Rothberger property.}
  \item{$\kappa < {\sf cov}(\mathcal{M}_{\reals})$.} 
\end{enumerate}
\end{proposition}
{\bf Proof:} $2\Rightarrow 1$: Let $X$ be a Lindel\"of space of cardinality $\kappa$. Let $(\mathcal{U}_n:n<\infty)$ be a sequence of open covers of $X$. Since $X$ is Lindel\"of we may assume each $\mathcal{U}_n$ is countable and enumerate it as $(U^n_k:k<\infty)$. For each $x$ in $X$ define
\[
  N_x = \{f\in\,^{\omega}\omega:(\forall n)(x\not\in U^n_{f(n)})\}.
\]
Each $N_x$ is nowhere dense in $^{\omega}\omega$ with the usual product topology, which is homeomorphic to the set of irrational numbers. Then $\bigcup_{x\in X}N_x\neq \,^{\omega}\omega$. Choose an $f\in\,^{\omega}\omega\setminus(\bigcup_{x\in X}N_x)$. Then the sequence $(U^n_{f(n)}:n<\infty)$ witnesses the Rothberger property of $X$ for $(\mathcal{U}_n:n<\infty)$.\\
$1\Rightarrow 2$: Note that 1 implies that every set of reals of cardinality $\kappa$ has the Rothberger property. Use the fact that ${\sf cov}(\mathcal{M}_{\reals})$ is the minimal cardinality of a set of reals that does not have the Rothberger property (\cite{FM}, Theorem 5).
$\epf$ 

Since ${\sf cov}(\mathcal{M}_{\reals})$ is the least $\kappa$ for which MA$_{\kappa}$(countable) is false, we have

\begin{corollary}\label{MAsmallLindelof} MA(countable) is equivalent to the statement that every Lindel\"of space of cardinality less than $2^{\aleph_0}$ has the Rothberger property.
\end{corollary}

\begin{center}{\bf Cardinals forbidden to support points ${\sf G}_{\delta}$ Rothberger topologies.}\end{center} 

Another application of ideas in the proof of Theorem \ref{rvmrothberger} gives: 
\begin{theorem}\label{wkcpt}
If Lebesgue measure can be extended to a countably additive measure on any collection of $2^{\aleph_0}$ subsets of the real line, then there is no points ${\sf G}_{\delta}$ Rothberger space of cardinality $2^{\aleph_0}$.
\end{theorem}
{\bf Proof:}
Suppose $X$ is a Lindel\"of space of cardinality $2^{\aleph_0}$ with all points ${\sf G}_{\delta}$. We may assume that as a set, $X$ is the closed unit interval. For each $x\in X$ choose a sequence $(U_n(x):n<\infty)$ of neighborhoods with $\{x\} = \bigcap_{n<\infty}U_n(x)$ and $U_{n+1}(x)\subseteq U_n(x)$, all $n$. Then $\{U_n(x):n<\infty,\, x\in X\}$ is a family of $2^{\aleph_0}$ subsets of $X$. Let $\mu$ be a countably additive extension of Lebesgue measure which also measures each $U_n(x)$.

For each $m$, choose for each $x$ an $n(m,x)$ such that $\mu(U_{n(m,x)}(x))< \frac{1}{2^{m+2}}$. Then
\[
  \mathcal{U}_m = \{U_{n(m,x)}(x):x\in X\}
\]
is an open cover of $X$. For the sequence $(\mathcal{U}_m:m<\infty)$ of open covers of $X$, consider any sequence $(U_m:m<\infty)$ with $U_m\in\mathcal{U}_m$ for each $m$. Then we have for each $m$ that $\mu(U_m)<\frac{1}{2^{m+2}}$, and so
\[
  \mu(\cup_{m<\infty}U_m)\le \sum_{m=1}^{\infty}\frac{1}{2^{m+2}} \le \frac{1}{2} < \mu(X).
\]
But then the sequence $(U_m:m<\infty)$ does not cover $X$. It follows that $X$ is not Rothberger.
$\epf$

T. Carlson proved that the consistency of the measure extension hypothesis of Theorem \ref{wkcpt} implies the consistency of the existence of a weakly compact cardinal, and Carlson and Prikry also derived the consistency of this measure extension hypothesis from  the consistency of the existence of a weakly compact cardinal. See Section 6 of \cite{TC} in this regard. Thus:

\begin{corollary}\label{wkcptrothb} If it is consistent that there is a weakly compact cardinal, then it is consistent that there is no points ${\sf G}_{\delta}$ Rothberger space of cardinality $2^{\aleph_0}$.
\end{corollary}

By the results of \cite{TC}, it is for example consistent that $2^{\aleph_0} = \aleph_2$, and for any family $\mathcal{F}$ of $\aleph_1$ subsets of the real line Lebesgue measure can be extended to a countably additive measure that measures also each element of $\mathcal{F}$. This, however, does not imply that there are no points ${\sf G}_{\delta}$ Rothberger spaces of size $\aleph_1$, as demonstrated by Section 4, Example 2. 

Carlson and Prikry derived the consistency of the measure extension hypothesis by starting with a model containing a weakly compact cardinal $\kappa$ and then adding $\kappa$ random reals. 
This leads to an alternative proof of Shelah's result \cite{SS} that no weakly compact cardinal can be topologized as a points ${\sf G}_{\delta}$ Lindel\"of space. First, observe:
\begin{corollary}\label{noweakcptrothberger}
There is no points ${\sf G}_{\delta}$ Rothberger topology on a weakly compact cardinal.
\end{corollary}
{\bf Proof:}
For suppose $X$ is a set of cardinality a weakly compact cardinal, and that $X$ carries a topology in which its points are ${\sf G}_{\delta}$, and it is a Rothberger space. Add $\vert X\vert$ random reals. By Theorem \ref{randomrothbergerpreserve} $X$ is still a Rothberger space with points ${\sf G}_{\delta}$ in the generic extension. But the generic extension is Carlson and Prikry's model for $\vert X\vert=2^{\aleph_0}$ and the extendibility of Lebesgue measure to any family of $2^{\aleph_0}$ subsets of the real line. Then Theorem \ref{wkcpt} gives the contradiction that $X$ is not a points ${\sf G}_{\delta}$ Rothberger space. 
$\epf$

\begin{corollary}[Shelah]\label{noweakcptlindelof}
There is no points ${\sf G}_{\delta}$ Lindel\"of topology on a weakly compact cardinal.
\end{corollary}
{\bf Proof:}
Suppose on the contrary that $\kappa$ is a weakly compact cardinal carrying a points ${\sf G}_{\delta}$ Lindel\"of topology.
Add $\aleph_1$ Cohen reals. The weakly compact cardinal $\kappa$ is still weakly compact, but by Theorem \ref{cohenrothberger} it now carries a points ${\sf G}_{\delta}$ Rothberger topology, contradicting Corollary \ref{noweakcptrothberger}.
$\epf$

Continuing the theme that we can prove for Rothberger spaces what we would like to prove for Lindel\"of spaces, we obtain in Theorem \ref{collapserothb} for Rothberger spaces results discussed but not achieved for Lindel\"of spaces in \cite{BT}.
For $\lambda<\kappa$ infinite regular cardinal numbers let $\mathbf{Lv}(\kappa,\lambda)$ be the partially ordered set whose elements are functions $p$ such that $\vert p\vert<\lambda$, $dom(p)\subseteq \kappa\times\lambda$ and for all $(\alpha,\xi)\in dom(p)$, $p(\alpha,\xi)\in\alpha$; elements $p$ and $q$ are ordered by $p \le q$ if, and only if, $q\subseteq p$. It is well-known that when $\kappa$ is strongly inaccessible, each antichain of ${\mathbf{Lv}}(\kappa,\lambda)$ is of cardinality less than $\kappa$. Also, as $\lambda$ is regular, ${\mathbf{Lv}}(\kappa,\lambda)$ is $\lambda$-closed. It is also well-known that in the generic extension obtained by forcing with ${\mathbf{Lv}}(\kappa,\lambda)$, $\kappa$ is a cardinal number, but is the successor of $\lambda$. In Theorem \ref{collapserothb} the phrase ``L\'evy-collapse $\cdots$ to $\omega_2$ with countable conditions" means ``Force with ${\mathbf{Lv}}(\kappa,\omega_1)$ where $\kappa$ is supercompact (measurable)."

The argument uses standard reflection methods (see e.g. \cite{DTW1}). Specific properties of supercompactness of a cardinal are at the core of these arguments, and are now recalled for the reader's convenience: A cardinal $\kappa$ is \emph{supercompact} if there is for each cardinal $\lambda\ge\kappa$ an elementary embedding $i:V\rightarrow M$ from the set-theoretic universe $V$ into a transitive class $M$ such that 
\begin{itemize}
  \item{for each $\xi<\kappa$, $i(\xi)=\xi$, but $i(\kappa)>\lambda$, and}
  \item{$^{\lambda}M \subseteq M$.}  
\end{itemize}
Here $^{\lambda}M$ denotes the class of sequences of length $\lambda$, where the terms of the sequences are elements of $M$. It follows by transfinite induction that for all $\alpha\le \lambda$ we have $M_{\alpha}=V_{\alpha}$ - that is, the cumulative hierarchy as computed in V and in M coincides at least up to $\lambda$.

\begin{theorem}\label{collapserothb} L\'evy-collapse a supercompact (measurable) cardinal to $\omega_2$ with countable conditions. Then every Rothberger (${\sf T}_2$) space of character $\le \aleph_1$ includes a Rothberger subspace of size $\le\aleph_1$.
\end{theorem}
{\bf Proof:} First, the supercompact version. Let $\kappa$ be a supercompact cardinal in the ground model. Let $G$ be ${\mathbf{Lv}}(\kappa,\omega_1)$-generic over $V$. In $V\lbrack G\rbrack$, let $X$ be a Rothberger space of character $\le \aleph_1$ and cardinality larger than $\aleph_1$. Let $\alpha$ be an initial ordinal such that $X$, all open covers of $X$, and all sequences of open covers of $X$, as well as $\mathbf{Lv}(\kappa,\omega_1)$ and any of its antichains are members of $V_{\alpha}\lbrack G\rbrack$. Choose a regular cardinal $\lambda> 2^{\alpha}$. In $V\lbrack G\rbrack$ let $\mu$ be the initial ordinal corresponding to $\vert X\vert$. Then $\mu$ is also a cardinal in $V$ and in $V$, $\kappa\le \mu < \alpha < \lambda$.

By supercompactness of $\kappa$ fix an elementary embedding $i:V\rightarrow M$ with $i(\kappa)>\lambda$ and $\,^{\lambda}M\subseteq M$. Since for each $\beta\le\lambda$ we have $M_{\beta}=V_{\beta}$ it follows that ${\mathbf{Lv}}(\kappa,\omega_1)$ and all $V$-antichains of it are elements of $M$. Thus $G\subseteq \mathbf{Lv}(\kappa,\omega_1)$ is $\mathbf{Lv}(\kappa,\omega_1)$-generic over $M$ if, and only if, it is over $V$. Moreover for all $\beta\le \lambda$, $V_{\beta}\lbrack G\rbrack = M_{\beta}\lbrack G\rbrack$. It follows that in fact $X\in M\lbrack G\rbrack$. Moreover, any family of open covers for $X$ in $V\lbrack G\rbrack$ is also in $M\lbrack G\rbrack$. Thus, $M\lbrack G\rbrack \models\ ``X \mbox{ is a Rothberger space of cardinality $>\aleph_1$ and character $\le \aleph_1$.}$" 

Since $\mathbf{Lv}(\kappa,\omega_1)$ has no antichain of cardinality $\kappa$, there is an $i(\mathbf{Lv}(\kappa,\omega_1))$-generic (over $M$) filter $G^*$ such that $p\in G$ implies $i(p)\in G^*$ (see Proposition 2.2 of \cite{DTW1}). But then $i$ extends to an elementary embedding $j:V\lbrack G\rbrack \rightarrow M\lbrack G^*\rbrack$ (see Proposition 2.1 of \cite{DTW1}).

Until further notice we now work in $M\lbrack G^*\rbrack$.  
The equation $i({\mathbf{Lv}}(\kappa,\omega_1)) = {\mathbf{Lv}}(i(\kappa),\omega_1) = {\mathbf{Lv}}(\kappa,\omega_1) \times {\mathbf{Lv}}(i(\kappa)\setminus\kappa,\omega_1)$ implies that $M\lbrack G^*\rbrack$ is of the form $M\lbrack G\rbrack\lbrack H\rbrack$ where $H$ is ${\mathbf{Lv}}(i(\kappa)\setminus\kappa,\omega_1)$-generic over $M\lbrack G\rbrack$. By Theorem \ref{countablyclosed} Rothberger is preserved by countably closed forcing, giving
  $M\lbrack G^*\rbrack\models ``X$ is a Rothberger space".

The bijection $j\lceil_X$ from $X$ to $j\lbrack X\rbrack$ induces a homeomorphic topology, say $\mathcal{T}$, on $j\lbrack X\rbrack$: Put a subset $U$ of $j\lbrack X\rbrack$ in $\mathcal{T}$ if its inverse image under $j$ is an open subset of $X$. Then $(j\lbrack X\rbrack,\mathcal{T})$ is a Rothberger space. 

The subset $j\lbrack X\rbrack$ of $j(X)$ inherits a topology from $j(X)$, say $\mathcal{S}$. Compare the two spaces $(j\lbrack X\rbrack,\mathcal{T})$ and  $(j\lbrack X\rbrack,\mathcal{S})$. First note that $\mathcal{S}$ contains sets of the form $j(U)\cap j\lbrack X\rbrack$, where $U\subseteq X$ is open in $X$. Conceivably, $\mathcal{S}$ may contain also other sets. Since $j(U)\cap j\lbrack X\rbrack = j\lbrack U\rbrack$ (see Lemma 2.0 of \cite{DTW1}) is an element of the topology $\mathcal{T}$ on $j\lbrack X\rbrack$, we have $\mathcal{T} \subseteq \mathcal{S}$. The character restriction ensures that $\mathcal{S} = \mathcal{T}$ (see the argument at the top of p. 45 of \cite{DTW1}). 
Thus, $j\lbrack X\rbrack$ is a Rothberger subspace of $j(X)$.

Since $j$ is an elementary embedding, $j(X)$ is a Rothberger space with character $\le \aleph_1$ and cardinality $j(\mu)$.
Since $j\lbrack X\rbrack$ is an uncountable subset of $j(X)$, we conclude that in $M\lbrack G^*\rbrack$ the statement 
\begin{quote}
   ``$j\lbrack X\rbrack \subseteq j(X)$ is an uncountable Rothberger subspace of $j(X)"$ 
\end{quote}
 as well as the statement $``j(\kappa) =\aleph_2 \mbox{ and } \vert\mu\vert=\aleph_1"$ are true. This implies: 
\begin{itemize}
  \item{$M\lbrack G^*\rbrack\models ``j(X)\mbox{ has a Rothberger subspace of cardinality } \aleph_1"$}
\end{itemize} 

This concludes working in $M\lbrack G^*\rbrack$. Since $j(\aleph_1) = \aleph_1$ and $j$ is an elementary embedding of $V\lbrack G\rbrack$ into $M\lbrack G^*\rbrack$, in $V\lbrack G\rbrack$ it is true that $X$ has a Rothberger subspace of cardinality $\aleph_1$. This concludes the proof for the supercompact case.

For the measurable version, note that a Lindelof ${\sf T}_2$ space of character $\le \aleph_1$ has cardinality $\le 2^{\aleph_1}$. All the cardinals below the measurable $\kappa$ will have their power sets collapsed, so $2^{\aleph_1} = \aleph_2$ - the former measurable - in the extension. Measurability yields a non-trivial elementary embedding moving $\kappa$, which again can be extended to a generic elementary embedding, whence we proceed as in the supercompact case.
$\epf$

In \cite{BvD} Baumgartner and van Douwen introduced the notion of a \emph{weakly measurable cardinal}: A cardinal number $\kappa$ is said to be weakly measurable if there is a sequence $(\mathcal{F}_n:n<\infty)$ of non-principal ultrafilters on $\kappa$ such that $\bigcap_{n<\infty}\mathcal{F}_n$ is a $\kappa$-complete filter on $\kappa$. It is known that the least such $\kappa$, if not measurable, is strictly less than $2^{\aleph_0}$ and that $\kappa \neq {\sf cof}(2^{\aleph_0})$. It is also known that the existence of a measurable cardinal is equiconsistent with the existence of a weakly measurable cardinal less than $2^{\aleph_0}$. The model given in Theorem 5.14 of \cite{BvD} for obtaining the consistency of a weakly measurable cardinal below $2^{\aleph_0}$ from the consistency of the existence of a measurable cardinal is as follows: Let $\kappa$ be a measurable cardinal in a model satisfying GCH. Let $\lambda>\kappa$ be a cardinal with cofinality not $\kappa$ and not $\omega$. Forcing with ${\mathbb P}(\lambda)$ (the partially ordered set for adding $\lambda$ Cohen reals) produces a model in which $\kappa<\lambda = 2^{\aleph_0}$ is weakly measurable. Thus, in some sense, weakly measurable cardinals are the Cohen analogue of real-valued measurable cardinals. In this model the set of Cohen reals is a Lusin set and thus has the Rothberger property. This shows that weakly measurable cardinals are not provably upper bounds on the possible cardinality of points ${\sf G}_{\delta}$ Rothberger spaces. Moreover, since every uncountable subset of a Lusin set is a Lusin set, every subset of this set of Cohen reals has the Rothberger property, and thus there are points ${\sf G}_{\delta}$ Rothberger spaces of all cardinalities less than or equal to $\lambda$ in this model.

\section{The Gerlits-Nagy property} 

Gerlits and Nagy \cite{GN} used the symbol $*$ to denote one of the several covering properties they introduced. We shall call $*$ the \emph{Gerlits-Nagy property}, and define it later below. Though the Gerlits-Nagy property is formally stronger than the Rothberger property the relationship between these two properties is somewhat complicated. Often ZFC examples of Rothberger spaces are also Gerlits-Nagy spaces. Some reasons for this phenomenon will be discussed below. An examination of Gerlits-Nagy spaces requires examining the Hurewicz- and Menger- properties which we introduce now.

The symbol $\sfin^{\alpha}(\mathcal{A},\mathcal{B})$ denotes the statement:
\begin{quote} For each sequence $(A_{\gamma}:\gamma<\alpha)$ of elements of $\mathcal{A}$, there is a sequence $(B_{\gamma}:\gamma<\alpha)$ of finite sets such that for each $\gamma$ we have $B_{\gamma}\subseteq A_{\gamma}$, and $\bigcup_{\gamma<\alpha}B_{\gamma}\in\mathcal{B}$. 
\end{quote}

The corresponding game, denoted $\gfin^{\alpha}(\mathcal{A},\mathcal{B})$, is played as follows: Players ONE and TWO play $\alpha$ innings. In the $\gamma$-th inning ONE chooses an $O_{\gamma}\in\mathcal{A}$, and TWO responds with a finite set $T_{\gamma}\subseteq O_{\gamma}$. A play
\[
  O_0,\, T_0,\, \cdots,\, O_{\gamma},\, T_{\gamma},\, \cdots
\]
is won by TWO if $\bigcup_{\gamma<\alpha}T_{\gamma}\in \mathcal{B}$; else, ONE wins.

$\sfin^{\omega}(\open,\open)$ is known as the \emph{Menger} property. W. Hurewicz \cite{H25} introduced this property and showed that in metrizable spaces it is equivalent to a basis property introduced by K. Menger \cite{Me}. Hurewicz also showed that a space has the property $\sfin^{\omega}(\open,\open)$ if, and only if, ONE has no winning strategy in the game $\gfin^{\omega}(\open,\open)$. An exposition of Hurewicz's result is given in Theorem 13 of \cite{coc1}. Clearly, the Rothberger property implies the Menger property.

Hurewicz \cite{H25} also introduced a property stronger than the Menger property, called the \emph{Hurewicz} property. The Hurewicz property is defined as follows: For each sequence $(\mathcal{U}_n:n<\omega)$ of open covers of $X$ there is a sequence $(\mathcal{V}_n:n<\omega)$ of finite sets such that for each $n$, $\mathcal{V}_n\subseteq\mathcal{U}_n$, and such that for each $x\in X$, for all but finitely many $n$, $x\in\bigcup\mathcal{V}_n$. The Hurewicz property can also be described as follows:

Call an open cover $\mathcal{U}$ of a space \emph{large} if each element of the space is contained in infinitely many members of $\mathcal{U}$. The symbol $\Lambda$ denotes the collection of large covers of a space. One can show that $\sfin^{\omega}(\open,\open)$ is equivalent to $\sfin^{\omega}(\Lambda,\Lambda)$.

A large cover $\mathcal{U}$ of a space $X$ is \emph{groupable} if there is a partition $\mathcal{U} = \bigcup_{n<\infty}\mathcal{U}_n$ such that each $\mathcal{U}_n$ is finite, for $m\neq n$ we have $\mathcal{U}_m\cap\mathcal{U}_n = \emptyset$, and for each $x\in X$, for all but finitely many $n$, $x\in \bigcup \mathcal{U}_n$.

By Theorem 12 of \cite{coc7} a Lindel\"of space has the Hurewicz property if, and only if, it has the Menger property and each countable large cover is groupable. 
$\sigma$-compactness implies the Hurewicz property, and the Hurewicz property implies the Menger property. Menger conjectured that in separable metric spaces the Menger property implies $\sigma$-compactness; Hurewicz conjectured that in separable metric spaces the Hurewicz property implies $\sigma$-compactness. It is known in ZFC that even in the class of separable metric spaces  none of the above implications is reversible: In \cite{CP} it is shown that the Menger property does not imply the Hurewicz property; in \cite{coc2} it is shown that the Hurewicz property does not imply $\sigma$-compactness. For more on these implications, see \cite{TZ}. Telg\'arsky \cite{Topsoe} proved for separable metric spaces: if TWO has a winning strategy in the game $\gfin^{\omega}(\open,\open)$, then that space is $\sigma$-compact (the converse is evidently true).

In Theorems 14 and 19 of \cite{NSW} it is shown that the Gerlits-Nagy property is equivalent to the Hurewicz property plus the Rothberger property. We take this characterization as our definition of the Gerlits-Nagy property.

\begin{center}{\bf Cohen reals and the Hurewicz property.}\end{center}

We have already seen that Cohen reals do not preserve not Rothberger. Thus, Cohen reals do not preserve not Menger. The situation for the Hurewicz property, and thus the Gerlits-Nagy property, is different. In the proof of Theorem \ref{nonhurewiczpreserve} we use the following well-known fact:
\begin{lemma}\label{cohenunbounded}
If $\mathcal{F}\subseteq\,^{\omega}\omega$ is unbounded in the eventual domination order and $\kappa>0$ is a cardinal number then
 \[
   {{\mathbf 1}}_{{\mathbb P}(\kappa)} \forces ``\check{\mathcal{F}} \mbox{ is unbounded }".
 \]
\end{lemma}
Though we have not found a direct reference (which undoubtedly exists), Lemma \ref{cohenunbounded} can be deduced from \cite{steprans} Lemma 3.1 together with the fact that for any ${\mathbb P}(\kappa)$-name $\dot{f}$ there is a countable set $S\subseteq\kappa$ with $\dot{f}$ in fact a ${\mathbb P}(S)$-name.

\begin{theorem}\label{nonhurewiczpreserve}
Let $\kappa>0$ be a cardinal number and let $X$ be a Lindel\"of space. If $X$ does not have the Hurewicz property then
 \[
   {{\mathbf 1}}_{{\mathbb P}(\kappa)} \forces ``\check{X} \mbox{ does not have the Hurewicz property }".
 \]
\end{theorem}
{\bf Proof:} 
 Let $X$ be a Lindel\"of space which is not Hurewicz. Let the sequence $(\mathcal{U}_n:n<\omega)$ of open covers of $X$ witness that $X$ is not Hurewicz. Since $X$ is Lindel\"of, we may assume each $\mathcal{U}_n$ is countable and has an enumeration $(U^n_k:k<\omega)$ such that for all $n$ and $k$, $U^n_k\subseteq U^n_{k+1}$. We shall show that 
\[
  {{\mathbf 1}}_{{\mathbb P}(\kappa)} \forces ``({\check{\mathcal{U}}_n}:n\in\check{\omega}) \mbox{ witnesses that }\check{X} \mbox{ is not Hurewicz}".
\]

Suppose the contrary and choose $p\in{\mathbb P}(\kappa)$ such that $p\forces ``\check{X} \mbox{ is Hurewicz}"$.
For each $x\in X$ define a function $f_x:\omega\longrightarrow\omega$ such that for each $n$, $f_x(n) = \min\{k:(\forall m\ge k)(x\in U^n_m)\}$. Note that the particular sequence of open covers of $X$ witness that $X$ is not Hurewicz if, and only if, the associated family of functions $\{f_x:\,x\in X\}$ is unbounded. Since $X$ is not Hurewicz, $\{f_x:x\in X\}$ is unbounded in $^{\omega}\omega$. 

However, $p\forces ``(\exists \tau\in{\,^{\omega}\omega})(\forall x\in \check{X})(f_x\prec \tau)"\footnote{$f\prec g$ denotes that for all but finitely many $n$, $f(n)<g(n)$.}$. Choose a ${\mathbb P}(\kappa)$ name $\dot{f}$ such that $p\forces ``(\forall x\in \check{X})(f_x\prec \dot{f})"$. Then choose a countable subset $C$ of $\kappa$ such that $\dot{f}$ is a ${\mathbb P}(C)$-name and $p\in {\mathbb P}(C)$. This gives a contradiction: ${\bf 1}_{{\mathbb P}(C)}\forces ``\{f_x:x\in \check{X}\} \mbox{ is unbounded}"$ (by Lemma \ref{cohenunbounded}) while $p\forces ``\{f_x:x\in \check{X}\} \mbox{ is bounded}"$.
$\epf$

The referee pointed out that the Cohen reals poset ${\mathbb P}(\kappa)$ does not preserve the Hurewicz property: This is most efficiently seen as follows: If $X$ is a co-dense subset of a complete metric space $M$ and if $X$ has the Hurewicz property, then $X$ is a meager subset of $M$ (see the proof of Theorem 5.5 of \cite{coc2}). But it is well-known that upon forcing with ${\mathbb P}(\omega)$, the set of ground model reals is a codense subset of the set of real numbers of the extension, but is not a meager subset of the reals of the extension (see Theorem 3.2 of \cite{steprans}). To now see this is the case for ${\mathbb P}(\kappa)$ also for uncountable $\kappa$, note: A ${\mathbb P}(\kappa)$-name for a comeager Borel set disjoint from the ground model reals is, for some countable subset $C$ of $\kappa$, a ${\mathbb P}(C)$-name.

There is a limited class of Hurewicz spaces for which we could prove that Hurewicz is preserved by Cohen real forcing. This class of Hurewicz spaces has appeared under different names in the study of small sets of reals. In \cite{pawl2} these are called ``property B", while in \cite{BJ} these are called ``sets in $\mathcal{H}$". Call a Hurewicz space $X$ \emph{strongly Hurewicz} if for each Borel function $F:X\rightarrow\,^{\omega}\omega$ there is a function $g\in\,^{\omega}\omega$ such that for each $x\in X$, for all but finitely many $n$, $F(x)(n)<g(n)$. Some properties of this class of Hurewicz spaces can be found in \cite{ST}, such as that such subsets of the real line are zero-dimensional and that Sierpi\'nski sets are strongly Hurewicz.

Recall the notion of an n-dowment. These are defined as follows: Fix a positive integer $n$. A family $\mathcal{L}$ of finite subsets of ${\sf Fn}(\kappa,2)$ is said to be an $n$-\emph{dowment} if for each maximal antichain $A\subseteq{\sf Fn}(\kappa,2)$ there is an $L\in\mathcal{L}$ with $L\subseteq A$, and for each $p\in{\sf Fn}(\kappa,2)$ with domain of cardinality $n$ and for any $L_1,\, L_2,\, \cdots,\, L_n \in\mathcal{L}$ there are $q_i\in L_i,\, i\le n$ such that $\{q_i:\, i\le n\}\cup\{p\}$ have a common extension in ${\sf Fn}(\kappa,2)$. For each $\kappa$ and each $n$ there is an $n$-dowment - see Lemma 1.1 of \cite{DTW1}.

Also recall: The Hurewicz game on a space $X$ is played as follows: ONE and TWO play an inning for each positive integer $n$. In the $n$-th inning ONE first chooses an open cover $O_n$ for $X$, and then TWO responds with a finite set $T_n\subseteq O_n$. A play $O_1,\, T_1,\, \cdots,\, O_n,\, T_n,\, \cdots$ is won by TWO if for each $x\in X$, for all but finitely many $n$, $x\in\bigcup T_n$; else, ONE wins. The game was introduced in \cite{coc1} where it was shown in Theorem 27 that $X$ is a Hurewicz space if, and only if, ONE has no winning strategy in the Hurewicz game on $X$.

\begin{theorem}\label{stronghurewiczpreserve} 
Let $\kappa>0$ be a cardinal number and let $X$ be a subspace of $\reals$. 
If $X$ is a strong Hurewicz space, then 
${{\mathbf 1}}_{{\mathbb P}(\kappa)} \forces ``\check{X} \mbox{ has the Hurewicz property }"$.
\end{theorem}
{\bf Proof:} Since a union of countably many Hurewicz spaces is a Hurewicz space, we may assume that $X$ is a subspace of $\lbrack 0,\, 1\rbrack$. For each positive integer $n$, let $L_n$ be an $n$-dowment, as defined above.
 
Since $X$ is a set of real numbers, it suffices to show that the following statement is true in generic extensions by ${\mathbb P}(\kappa)$:
\begin{quote} For each ${\sf G}_{\delta}$-set $G\subseteq\lbrack 0,\, 1\rbrack$ with $X\subseteq G$ there is an ${\sf F}_{\sigma}$-set $F$ such that $X\subseteq F\subseteq G$. (See Theorem 5.7 of \cite{coc2}\footnote{This characterization of Hurewicz, suitably reformulated, holds in greater generality: \cite{BZ} extended the result to separable metric spaces, and \cite{FDT2} extended it to Lindel\"of ${\sf T}_3$ spaces.}.)
\end{quote}
Thus, let a sequence $(\dot{G}_n:n<\omega)$ of ${\mathbb P}(\kappa)$-names be given such that 
\[
{{\mathbf 1}}_{{\mathbb P}(\kappa)} \forces ``(\forall n)(\dot{G}_n\subseteq\lbrack0,\, 1\rbrack \mbox{ is open and includes }\check{X} \mbox{ and } \dot{G}_n\supseteq \dot{G}_{n+1})"
\]
Then for each $n$ let $\dot{C}_n$ be a ${\mathbb P}(\kappa)$-name such that  
\[
  {{\mathbf 1}}_{{\mathbb P}(\kappa)} \forces ``(\forall n)(\dot{C}_n = \lbrack0,\, 1\rbrack\setminus\dot{G}_n \mbox{ is compact and }\check{X}\,\cap\,\dot{C}_n =\emptyset)"
\]
Then for each $n$, and for each $x\in X$,
\[
  {{\mathbf 1}}_{{\mathbb P}(\kappa)} \forces ``(\exists\dot{U}_n)(\exists \dot{V}_n)(\dot{U}_n \mbox{ and }\dot{V}_n \mbox{ open, and } \overline{\dot{U}_n}\,\cap\,\overline{\dot{V}_n} = \emptyset \mbox{ and }x\in \dot{U}_n,\, \dot{C}_n\subseteq\dot{V}_n)"
\]
Choose for each $x\in X$ and each $n<\omega$ a maximal antichain $A(x,\dot{C}_n)\subseteq {\mathbb P}(\kappa)$, and for each $p\in A(x,\dot{C}_n)$ a clopen (in $X$) interval $U_p(x,\dot{C}_n)$, a ${\mathbb P}(\kappa)$-name $\dot{V}_{p,x}(\dot{C}_n)$ and a finite set $Q_{p,x}(\dot{C}_n)$ of pairs $r<s$ of rational numbers such that with $\dot{I}(r,s)$ the ${\mathbb P}(\kappa)$-name for the open interval $(\check{r},\check{s})$, $p$ forces the conjunction of the following three statements:
\begin{enumerate}
  \item{$``\dot{V}_{p,x}(\dot{C}_n) = \bigcup\{\dot{I}(r,s):\,{(r,s)\in {Q_{p,x}(\dot{C}_n)}}\}"$,}
  \item{$``\overline{\check{U}_p(x,\dot{C}_n)}\, \cap\, \overline{\dot{V}_{p,x}(\dot{C}_n)} = \emptyset"$,}
  \item{$``\dot{C}_n \subseteq \dot{V}_{p,x}(\dot{C}_n)"$.}
\end{enumerate} 
In the ground model it is true that $U_p(x,\dot{C}_n)$ is disjoint from $\cup\{(r,s):(r,s)\in Q_{p,x}(\dot{C}_n)\}$, because ${\mathbf 1}_{{\mathbb P}(\kappa)}\forces``(r,s)\check{}\subseteq \dot{I}(r,s)"$.

Then for each $x,\, n \mbox{ and }k$ choose a finite set $F_k(x,\dot{C}_n)\subseteq A(x,\dot{C}_n)$ which is a member of the $k$-dowment $L_k$. For fixed $n$ and $k$ define for each $x$
\begin{center}
  \begin{tabular}{l}
   $U^n_k(x) = \cap\{U_p(x,\dot{C}_n):p\in F_k(x,\dot{C}_n)\}$\\
   $\dot{V}^n_{k,x} = \cap\{\dot{V}_{p,x}(\dot{C}_n):p\in F_k(x,\dot{C}_n)\}$\\
  \end{tabular}
\end{center}
In the second case we have defined a ${\mathbb P}(\kappa)$-name. For there is in the ground model the corresponding finite set $Q^n_{k,x}$ of pairs $(r,s)$ of rational numbers obtained from the corresponding finite sets $Q_{p,x}(\dot{C}_n)$, $p\in F_k(x,\dot{C}_n)$ (by taking the set of intersections of the finitely many rational intervals associated with the sets $Q_{p,x}(\dot{C}_n)$) such that $\dot{V}^n_{k,x} = \cup\{\dot{I}(r,s):(r,s)\in Q^n_{k,x}\}$. Define in the ground model the set $v^n_{k,x} = \cup_{(r,s)\in Q^n_{k,x}}(r,s)$. 
Observe that for each $p\in F_k(x,\dot{C}_n)$ we have $p\forces ``\overline{\check{U}^n_k(x)}\,\,\cap\,\, \overline{\dot{V}^n_{k,x}} = \emptyset"$, and in the ground model we have $\overline{U^n_k(x)}\,\, \cap\,\, \overline{v^n_{k,x}} = \emptyset$.

Next, for each $n$ and $k$ define $\mathcal{U}_{n,k} = \{U^n_k(x):x\in X\}$, a ground-model clopen cover of $X$. For each fixed $n$ define a strategy $\sigma_n$ of player ONE of the Hurewicz game on $X$ as follows:
{\flushleft{\bf Definition of }$\sigma_n(\emptyset)$:} $\sigma_n(\emptyset) = \mathcal{U}_{n,1}$.\\
Suppose that TWO responds with the finite set $T^n_1\subseteq \sigma_n(\emptyset)$. 
{\flushleft{\bf Definition of }$\sigma_n(T^n_1)$:} Fix a finite set $\{x^{n,1}_1,\cdots,x^{n,1}_{\ell^n_1}\}\subseteq X$ with $T^n_1 = \{U^n_1(x^{n,1}_i):i\le \ell^n_1\}$. Define $m^n_1 = \max\{\vert p\vert: p\in \cup_{i\le \ell^n_1}F_1(x^{n,1}_i,\dot{C}_n)\}$. Then put
\[
  \sigma_n(T^n_1)=\mathcal{U}_{n,m^{n}_1+1}.
\]
Suppose that TWO responds with the finite set $T^n_2\subseteq \sigma_n(T^n_1)$. 
{\flushleft{\bf Definition of }$\sigma_n(T^n_1, T^n_2)$:} Fix a finite set $\{x^{n,2}_1,\cdots,x^{n,2}_{\ell^n_2}\}\subseteq X$ with $T^n_2 = \{U^n_{m_1+1}(x^{n,2}_i):i\le \ell^n_2\}$. Define $m^n_2 = \max\{\vert p\vert: p\in \cup_{i\le \ell^n_2}F_{m^n_1+1}(x^{n,2}_i,\dot{C}_n)\}$. Then put
\[
  \sigma_n(T^n_1,T^n_2)=\mathcal{U}_{n,m^{n}_1+m^n_2+1},
\]
and so on.

Since $X$ has the Hurewicz property, $\sigma_n$ is not a winning strategy for ONE. Choose a $\sigma_n$-play lost by ONE, and let TWO's sequence of moves in this play be 
\[
  T^n_1,\, T^n_2,\, \cdots,\, T^n_k,\, \cdots
\]
Then we have for each $x\in X$ that for all but finitely many $k$, $x\in\cup T^n_k$. For this sequence of moves by TWO fix the symbols $x^{n,k}_i$, $\ell^n_k$ and $m^n_k$ as above in the definition of the strategy $\sigma_n$.

Note that for each $n$ we have: For each $j$, if $p\in F_{m^n_1+\cdots+m^n_{j-1}+1}(x^{n,j}_i,\dot{C}_n)$ for some $i\le \ell^n_j$, then $p\forces `` \overline{\cup\check{T}^n_j}\,\,\cap\,\, \overline{\dot{V}^n_j} = \emptyset"$, where henceforth $\dot{V}^n_j$ is a ${\mathbb P}(\kappa)$-name for $\dot{V}^n_{j-1}\cap(\cap_{i\le \ell^n_j}\dot{V}^n_{j,x^{n,j}_i})$.

Now consider an arbitrary $p\in{\mathbb P}(\kappa)$, and an arbitrary $n$. Consider any $k$ so large that $\vert p\vert < m^n_1+\cdots+ m^n_{k-1}+1$. Then for each $x^{n,k}_i$, $i\le \ell^n_k$, find a $q_i\in F_{m^n_1+\cdots+m^n_{k-1}+1}(x^{n,k}_i,\dot{C}_n)$ such that each pair $\{q_i,\, p\}$ has a common extension, say $r_i$. Then 
$r_i\forces ``(\forall t\ge k)(\overline{\cup\check{T}^n_k}\,\,\cap\,\,\overline{\dot{V}^n_t} = \emptyset)."$
It follows that 
\[
  p\forces``(\exists S_n\in\lbrack\check{\omega}\rbrack^{\aleph_0})(\forall k\in S_n)(\forall t\ge k)(\overline{\cup\check{T}^n_k}\,\,\cap\,\,\overline{\dot{V}^n_t} = \emptyset)"
\]
Since $p$ and $n$ were arbitrary we find that in fact
\[
  {\mathbf 1}_{{\mathbb P}(\kappa)}\forces``(\forall n)(\exists S_n\in\lbrack\check{\omega}\rbrack^{\aleph_0})(\forall k\in S_n)(\forall t\ge k)(\overline{\cup\check{T}^n_k}\,\,\cap\,\,\overline{\dot{V}^n_t} = \emptyset)"
\]

In the ground model define a function $F:X\rightarrow\,^{\omega}\omega$ as follows:
\[
  F(x)(n)=\min\{k:(\forall j\ge k)(x\in\cup T^n_j)\}.
\]
Since for each basic open subset $\lbrack(n_1,\cdots,n_j)\rbrack$ of $^{\omega}\omega$ we have
\[
  F^{-1}\lbrack\lbrack(n_1,\cdots,n_j)\rbrack\rbrack = \cap_{i\le j}((\cap_{t\ge n_i}(\cup T^i_j))\setminus (\cup T^i_{n_i-1})), 
\]
a Borel subset of $X$, and thus $F$ is a Borel function. Since $X$ is strongly Hurewicz we fix an increasing $g\in\,^{\omega}\omega$ such that for each $x\in X$, for all but finitely many $n$, $F(x)(n) < g(n)$.

Let $G$ be ${\mathbb P}(\kappa)$-generic. In $V\lbrack G\rbrack$ we do the following:

For each $n$ choose an infinite subset $S_n$ of $\omega$ such that for each $k\in S_n$, for all $t\ge k$, 
$\overline{\cup{T}^n_k}\,\,\cap\,\,\overline{{V}^n_t} = \emptyset$. Then, for each $n$, choose $f(n)\in S_n$ with $f(n) > g(n)$. 

It follows that
\begin{itemize} 
  \item{For each $x$, for all but finitely many $n$, $x\in \cup T^n_{f(n)}$.}
  \item{For each $n$, for all $t\ge f(n)$, $\overline{\cup{T}^n_{f(n)}}\,\,\cap\,\,\overline{{V}^n_t} = \emptyset.$}
\end{itemize}
For each $k$ put $F_k = \cap_{n\ge k}\overline{\cup{T}^n_{f(n)}}$, a closed set. By the second item above $F_k$ is disjoint from $\cup_{n\ge k}C_n$, and thus included in $\cap_{n\ge k}G_n$. Since each $G_k$ is open it is an ${\sf F}_{\sigma}$-set and thus, setting $E_k = F_k \cap G_1\cap\cdots\cap G_{k-1}$, we find that $E_k\subseteq G$ is an ${\sf F}_{\sigma}$-set. Also observe that $X\cap E_k = X\cap F_k$ since $X\subseteq G$.

By the first item above, $X\subseteq\cup_{k<\omega}E_k$, where the latter is an $F_{\sigma}$-set. It follows that in $V\lbrack G\rbrack$ the statement ``there is an $F_{\sigma}$-set $E$ with $X\subseteq E\subseteq G$" is true. Since $G$ was an arbitrary ${\mathbb P}(\kappa)$-generic filter, we have now proven
\[
{\mathbf 1}_{{\mathbb P}(\kappa)}\forces ``\check{X} \mbox{ is a Hurewicz space in the topology inherited from }\dot{\reals}."
\]
But since in the generic extension the topology of $X$ generated by the ground model topology of $X$ is a subset of the topology $X$ inherits from the real line of the generic extension, it finally follows that\footnote{In general, if $\tau_1\subseteq \tau_2$ are topologies on a space $X$ and $(X,\tau_2)$ has the Lindel\"of, Rothberger, Hurewicz or Menger property, then so does $(X,\tau_1)$.}
\[
{\mathbf 1}_{{\mathbb P}(\kappa)}\forces ``\check{X} \mbox{ is a Hurewicz space". }\epf
\]

\begin{problem}\label{stronghurewicz} Does Cohen real forcing preserve the strong Hurewicz property?
\end{problem}

Adding enough Cohen reals converts ground-model strongly Hurewicz spaces to Gerlits-Nagy spaces in the generic extension:
\begin{corollary}\label{Hurewicztogerlitsnagy}
Let $\kappa$ be an uncountable cardinal. For a space $X\subseteq\reals$,  if $X$ has the strong Hurewicz property, then
${{\mathbf 1}}_{{\mathbb P}(\kappa)} \forces ``\check{X} \mbox{ has the Gerlits-Nagy property }"$.
\end{corollary}
{\bf Proof:} The Hurewicz property implies the Lindel\"of property. Now apply Theorem \ref{cohenrothberger} and Theorem \ref{stronghurewiczpreserve} to conclude that a ground-model strongly Hurewicz subspace of $\reals$ is converted to a Rothberger space which still has the Hurewicz property.
$\epf$

\begin{center}{\bf Random reals and the Hurewicz and Menger properties}\end{center}

\begin{theorem}\label{nonmenger}
Let $\kappa$ be an uncountable cardinal number. If $X$ is a Lindel\"of space, the following are equivalent:
\begin{enumerate}
  \item{$X$ has the Menger property.} 
  \item{${{\mathbf 1}}_{{\mathbb B}(\kappa)} \forces ``{\check{X}} \mbox{ has the Menger property}".$}
\end{enumerate}
\end{theorem}
{\bf Proof:} Suppose that $X$ is a Lindel\"of space. \\
{\bf 2$\Rightarrow $1:} Suppose that for some positive element $b$ of ${\mathbb B}(\kappa)$,
\[
  b \forces ``\check{X} \mbox{ has the Menger property }\sfin^{\omega}(\open,\open)".
\]
 Let $(\mathcal{U}_n:n<\omega)$ be a ground-model sequence of open covers of $X$. We may assume each $\mathcal{U}_n$ is countable and enumerate it as $(U^n_m:m<\omega)$. Then 
\[
  b \forces ``(\check{\mathcal{U}}_n:n<\omega) \mbox{ is a sequence of open covers of } \check{X}"
\]
 and thus
\[
  b \forces ``(\exists \phi\in\,^{\omega}\omega)(\forall x \in \check{X})(\exists^{\infty}_n)(x\in \bigcup_{j\le\phi(n)}U^n_j)"
\]

Now it is well-known that there is an $f\in\, ^{\omega}\omega$ such that 
\[
  b \forces ``(\forall^{\infty}_n)(\phi(n) < f(n))".
\]
Fix such an $f$: Then
\[
   b \forces ``(\forall x \in \check{X})(\exists^{\infty}_n)(x\in \bigcup_{j\le f(n)}U^n_j)"
\]
But all parameters in this last sentence are from the ground model, thus the sentence is true in the ground model, and we find an $f\in\,^{\omega}\omega$ such that setting $\mathcal{V}_n = \{U^n_j:j\le f(n)\}$, the sequence of $\mathcal{V}_n$'s cover $X$. We conclude that $X$ has the Menger property in the ground model.\\
{\bf 1$\Rightarrow $2:} Consider a ${{\mathbb B}(\kappa)}$-name $(\mathcal{U}_n:n<\infty)$ such that
\[
  {{\mathbf 1}}_{{\mathbb B}(\kappa)}\forces``(\mathcal{U}_n:n<\infty) \mbox{ is a sequence of open covers of }\check{X}"
\]
Since the ground-model topology is a basis for the topology of X in the generic extension, we may assume that 
\[
  {{\mathbf 1}}_{{\mathbb B}(\kappa)} \forces ``\mbox{Each }\mathcal{U}_n\mbox{ consists of sets open in the ground-model topology}"
\]
For each $x\in X$ and each $n$, choose a maximal antichain $A^n_x\subseteq {\mathbb B}(\kappa)$ such that whenever $b\in A^n_x$ then there is a neighborhood $V_b^n(x)$ of $x$ such that $b\forces``(\exists U\in\mathcal{U}_n)(V_b^n(x)\subseteq U)"$.
Then choose a finite set $L^n_x\subseteq A^n_x$ such that $\mu(\bigcup L^n_x)\ge 1 - (\frac{1}{2})^n$. Then put
\[
  W_n(x) = \bigcap\{V^n_b(x): b\in L^n_x\}.
\]
The set $\mathcal{G}_n:=\{W_n(x):x\in X\}$ is a ground-model open cover of $X$. \\
{\bf Claim 1:} $(\forall x\in X)(\forall b\in {\mathbb B}(\kappa))(\exists m<\infty)(\forall n\ge m)((\exists V\in\mathcal{G}_n)(x\in V)\Rightarrow (\exists r\le b)(r\forces``(\exists U\in\mathcal{U}_n)(V\subseteq U)")).$

To see this, consider an $x\in X$ and a $b\in{\mathbb B}(\kappa)$. Choose $m$ so large that $\frac{1}{2^m} < \mu(b)$. Consider any $n\ge m$. There are $V\in\mathcal{G}_n$ with $x\in V$. Fix such a $V$, say $V = W_n(y)$, some $y\in X$. Then we have:
\[
  \mu(\Vert(\exists U\in\mathcal{U}_n)(V\subseteq U)\Vert)) \ge 1 - (\frac{1}{2})^n.
\]
Set $r = b \cap \Vert(\exists U\in\mathcal{U}_n)(V\subseteq U)\Vert$. This completes the proof of Claim 1.

Observe that the ``m" in Claim 1 can be replaced with any larger value. 
Since $X$ has the Menger property in the ground model, choose there for each $n$ a finite set $\mathcal{H}_n\subseteq\mathcal{G}_n$ such that
\[
  (\forall x\in X)(\exists^{\infty}_n)(x\in\bigcup\mathcal{H}_n).
\]
{\flushleft{\bf Claim 2:} $D^m_x = \{b\in {\mathbb B}(\kappa): (\exists n >m)(\exists V \in \mathcal{H}_n)(x\in V \mbox{ and } b\forces ``(\exists U \in {\check{{\mathcal{U}_n}}})(V\subseteq U)")\}$ is dense in ${\mathbb B}(\kappa)$.}

To see this, consider $x\in X$ and $m<\infty$. Consider any $b\in{\mathbb B}(\kappa)$. Choose by Claim 1 an $M>m$ such that for any $n\ge M$, if there is a $V\in\mathcal{G}_n$ with $x\in V$, then there is an $r\le b$ such that $r\forces``(\exists U \in \mathcal{U}_n)(V\subseteq U)"$. Then choose $n>M$ with $x\in\mathcal{H}_n$. Apply Claim 1 to obtain a member $q$ of $D^m_x$ with $q<b$. This completes the proof of Claim 2.

Now let $G$ be a ${\mathbb B}(\kappa)$-generic filter. For each $m$ and $x\in X$ choose a $b(m,x)\in D^m_x\cap G$. Then choose $n(m,x)>m$ and $H(m,x)\in\mathcal{H}_{n(m,x)}$ with $x\in H(m,x)$. Now $b(m,x)\forces ``(\exists U\in\mathcal{U}_{n(m,x)})(H(m,x)\subseteq U)"$, and thus this forced sentence is true in the generic extension by $G$. In the generic extension choose for each such $n$, and for each $V\in\mathcal{H}_n$ for which it is possible, a $U_V\in\mathcal{U}_n$ with $V\subseteq U_V$, and let $\mathcal{J}_n$ be the set of such chosen $U_V$'s. Then in the generic extension: For each $n$, $\mathcal{J}_n$ is a finite subset of $\mathcal{U}_n$, and $\bigcup_{n<\infty}\mathcal{J}_n$ is an open cover of $X$. It follows that in the generic extension $X$ has the Menger property. Since this holds for all ${\mathbb B}(\kappa)$-generic filters, 1$\Rightarrow$ 2 is proven.
$\epf$

Next we take up preservation and non-preservation of the Hurewicz property under forcing with ${\mathbb B}(\kappa)$.
\begin{theorem}\label{randomnonhurewicz} For $X$ a topological space: If ${{\mathbf 1}}_{{\mathbb B}(\kappa)} \forces ``\check{X} \mbox{ has the Hurewicz property}"$, then $X$ has the Hurewicz property
\end{theorem}
{\flushleft{\bf Proof:}} The proof of Theorem \ref{nonhurewiczpreserve} adapts to this case.
$\epf$

Unlike the case of Cohen reals forcing, for ${\sf T}_1$ regular Lindel\"of spaces the Hurewicz property is preserved by random real forcing. In preparation for proving this, note first that regular Lindel\"of spaces are normal, thus completely regular. Since such a space is ${\sf T}_1$, it embeds as a subspace of $\lbrack 0,\, 1\rbrack^{\kappa}$ for some large enough cardinal $\kappa$. Since homeomorphisms preserve both the Hurewicz and the non-Hurewicz properties, we may restrict attention to subsets of $\lbrack 0,\, 1\rbrack^{\kappa}$. First, we prove the result for $\kappa=\omega$, the case of separable metric spaces. Then we return to the more general case. For the special case of a separable metric space we use the following alternate characterization of the Hurewicz property of subsets of $\lbrack 0,\, 1\rbrack^{\naturals}$, the Hilbert cube ${\mathbb H}$:
\begin{proposition}[\cite{coc2}]\label{hurewiczgdelta} A subspace $X$ of $\hilbertcube$ has the Hurewicz property if, and only if, for each ${\sf G}_{\delta}$ subset $G$ of $\hilbertcube$ with $X\subseteq G$, there is an ${\sf F}_{\sigma}$ subset $F$ of $\hilbertcube$ with $X\subseteq F\subseteq G$.
\end{proposition}
{\flushleft{\bf Proof:}} See Theorem 5.7 of \cite{coc2}: The proof given there is for the real line, but adapts easily to the Hilbert cube. $\epf$

More generally this argument also applies to characterizing Hurewicz subsets of $\lbrack 0,\, 1\rbrack^{\lambda}$. We shall also need the following observation:
\begin{lemma}\label{gdeltasets} Let $\lambda$ be an uncountable cardinal number and let $X$ be a Lindel\"of subset of $\lbrack 0,\, 1\rbrack^{\lambda}$. For each ${\sf G}_{\delta}$-set $G\subseteq\lbrack 0,\, 1\rbrack^{\lambda}$ such that $X\subseteq G$ there is a ${\sf G}_{\delta}$ set $H$ and a countable set $C\subseteq \lambda$ such that 
\begin{enumerate}
  \item{$X\subseteq H\subseteq G$ and}
  \item{$H$ is homeomorphic to $\Pi_C\lbrack H\rbrack\times \lbrack 0,\, 1\rbrack^{\lambda\setminus C}$.}
\end{enumerate}
\end{lemma}
{\bf Proof:} First, consider an open set $U\subseteq\lbrack 0,\, 1\rbrack^{\lambda}$ with $X\subseteq U$. Choose for each $x\in X$ a basic open set $B(U,x)$ of the form $\Pi_{\alpha\in \lambda}S_{\alpha}\subseteq U$ where 
\[
  F(U,x) = \{\alpha\in\lambda:S_{\alpha}\mbox{ is a proper open subset of $\lbrack 0,\, 1\rbrack$ with rational endpoints }\}
\]
 is finite. Then $\{B(U,x):x\in X\}$ is an open cover of $X$ and thus has a countable subset, $\{B(U,x_n):n<\omega\}$ which covers $X$. Then $C_U = \cup_{n<\omega}F(U,x_n)$ is a countable subset of $\lambda$, and $V = \cup_{n<\omega}B(U,x_n)$ is an open set with $X\subseteq V\subseteq U$, and for any countable set $D$ with $C_U\subseteq D\subseteq \lambda$, $V$ is homeomorphic to $\Pi_{D}\lbrack V\rbrack \times \lbrack 0,\, 1\rbrack^{\lambda\setminus D}$.

Next, let $G$ be a ${\sf G}_{\delta}$ set containing $X$, and choose open sets $U_1\supseteq U_2\supseteq\cdots\supseteq U_n\supseteq \cdots$ such that $G = \cap_{n<\omega}U_n$. Applying the preceding remarks consecutively to each $U_n$ we find open sets $V_n$ and countable sets $C_n\subseteq \lambda$ such that for each $n$ we have
\begin{enumerate}
  \item{$U_n\supseteq V_n\supseteq V_{n+1}\supseteq X$;}
  \item{$C_n\subseteq C_{n+1}$;}
  \item{For each countable $D$ with $C_n\subseteq D\subseteq \lambda$, $V_n$ is homeomorphic to $\Pi_{D}\lbrack V_n\rbrack \times \lbrack 0,\, 1\rbrack^{\lambda\setminus D}$.} 
\end{enumerate}
Finally put $H=\cap_{n<\omega}V_n$, and $C = \cup_{n<\omega}C_n$.
$\epf$

\begin{theorem}\label{nonhurewicz}
Let $\kappa$ be an uncountable cardinal number. If $X$ is a {\sf T}$_3$ Hurewicz space then 
\[
 {{\mathbf 1}}_{{\mathbb B}(\kappa)} \forces ``\check{X} \mbox{ has the Hurewicz property}".
\]
\end{theorem}
{\flushleft{\bf Proof:}} First we treat the case of subspaces of the Hilbert cube ${\mathbb H}$. Thus, assume $X\subseteq \hilbertcube$ has the Hurewicz property. Assume that
\[
  {\mathbf 1}_{{\mathbb B}(\kappa)}\forces ``\dot{G}\subseteq\dot{\hilbertcube} \mbox{ is a {\sf G}$_{\delta}$-set and }\check{X}\subseteq \dot{G}."
\]
Fix ${\mathbb B}(\kappa)$-names $\dot{G}_n$, $n<\omega$, such that
\[
  {\mathbf 1}_{{\mathbb B}(\kappa)}\forces ``(\forall n)(\dot{G}_n\supseteq \dot{G}_{n+1}\supseteq \dot{G} \mbox{ are open and }\dot{G} = \bigcap_{n<\omega} \dot{G}_n."
\]
We may assume that for each $n$ we have ${\mathbb B}(\kappa)$-names $\dot{C}_n$ such that:
\[
  {\mathbf 1}_{{\mathbb B}(\kappa)}\forces ``(\forall n)(\dot{C}_n = \dot{\hilbertcube}\setminus \dot{G}_{n}, \mbox{ a nonempty compact set})."
\]

Let an $\epsilon>0$ be given. Consider a fixed $n$. For each $x\in X$ we have ${\mathbf 1}_{{\mathbb B}(\kappa)}\forces``\check{x}\not\in\dot{C}_n"$. For each $x$ choose a basic neighborhood $U_n(x)$ and a ${\mathbb B}(\kappa)$-name $\dot{V}_{n,x}$ so that
\[
  \mu(\Vert(\dot{C}_n\subseteq \dot{V}_{n,\check{x}} \mbox{ open, and }\check{U}_n(\check{x})\cap\dot{V}_{n,\check{x}} = \emptyset)\Vert) \ge 1-\frac{\epsilon}{2^n}.
\]
Now $\mathcal{U}_n = \{U_n(x):x\in X\}$ is an open cover of $X$.

Applying the fact that $X$ is Hurewicz, choose for each $n$ a finite set $\mathcal{F}_n\subseteq\mathcal{U}_n$ such that for each $x\in X$, for all but finitely many $n$, $x\in\bigcup\mathcal{F}_n$. For each $x\in X$ choose an $N(x)\in\naturals$ such that for all $n\ge N(x)$ we have $x\in\bigcup\mathcal{F}_n$.

For each $n$ let $I_n$ be a finite set for which $(U_n(x_i):i\in I_n)$ is a bijective enumeration of $\mathcal{F}_n$. Then for each $n$ define a ${\mathbb B}(\kappa)$-name $\dot{V}_n$ so that ${\mathbf 1}_{{\mathbb B}(\kappa)}\forces ``\dot{V}_n = \bigcap_{i\in \check{I}_n}\dot{V}_{n,\check{x}_i}$. For each $n$ we have $\mu(\Vert (\bigcup\check{\mathcal{F}}_n) \bigcap \dot{V}_n = \emptyset\Vert)\ge 1 - \frac{\epsilon}{2^n}.$ 

For fixed $N$ define $X_N=\{x\in X: N(x) = N\}.$ Then, for all $n\ge N$, $X_N\subseteq \bigcup\mathcal{F}_n$, and so, for each $n\ge N$, $\mu(\Vert \check{X}_N\cap \dot{V}_n = \emptyset\Vert) \ge 1 - \frac{\epsilon}{2^n}$ and consequently $\mu(\Vert \check{X}_N\cap (\bigcup_{n\ge N}\dot{V}_n) = \emptyset\Vert) \ge 1 - \epsilon\sum_{n=N}^{\infty}\frac{\epsilon}{2^n} = 1 - \frac{\epsilon}{2^{N-1}}.$ This gives $\mu(\Vert \check{X}_N\subseteq \dot{\hilbertcube}\setminus(\bigcup_{n\ge N}\dot{V}_n)\Vert) \ge 1 - \frac{\epsilon}{2^{N-1}}.$

Consequently we find:  
\[
  \mu(\Vert \check{X}_N \subseteq (\dot{\hilbertcube}\setminus(\bigcup_{n\ge N}\dot{V}_n)\cap (\bigcap_{n<N}\dot{G}_n)) \subseteq \dot{G} \Vert) \ge 1 - \frac{\epsilon}{2^{N-1}}.
\]
But then, since $\Vert\dot{\hilbertcube}\setminus(\bigcup_{n\ge N}\dot{V}_n) \mbox{ is $\sigma$-compact and }\bigcap_{n<N}\dot{G}_n \mbox{ is $\sigma$-compact}\Vert = {\mathbf 1}_{{\mathbb B}(\kappa)}$, we find that there is a ${\mathbb B}(\kappa)$-name $\dot{F}_N$ such that
\[
 \mu(\Vert(\dot{F}_N \mbox{ is $\sigma$-compact and }\check{X}_N\subseteq \dot{F}_N\subseteq \dot{G})\Vert \ge 1 - \frac{\epsilon}{2^{N-1}}.
\] 
Since we also have $\mu(\Vert \check{X}\subseteq \bigcup_{n=1}^{\infty}\dot{F}_n\Vert)\ge 1 - \epsilon\cdot(\sum_{n=1}^{\infty}\frac{1}{2^{n-1}}) = 1 - 2\cdot\epsilon$, we conclude that $\mu(\Vert(\exists\mbox{ a $\sigma$-compact set }\dot{F})(\check{X}\subseteq \dot{F}\subseteq \dot{G}\Vert)\ge  1 - 2\cdot\epsilon.$
As $\epsilon>0$ was arbitrary it follows that 
\[
{\mathbf 1}_{{\mathbb B}(\kappa)}\forces ``(\forall\mbox{ G$_{\delta}$ set }\dot{G} \subseteq \dot{\hilbertcube})(\exists \mbox{ an ${\sf F}_{\sigma}$-set }\dot{F})(\check{X}\subseteq \dot{F}\subseteq \dot{G})".
\]
By Proposition \ref{hurewiczgdelta} we have ${\mathbf 1}_{{\mathbb B}(\kappa)}\forces ``\check{X}\mbox{ is Hurewicz}".$

Next we consider the general case: Let $X$ be a ${\sf T}_3$ Hurewicz space and let $\lambda$ be the least infinite cardinal such that $X$ is homeomorphic to a subspace of $\lbrack 0,\, 1\rbrack^{\lambda}$. We may assume that $\lambda$ is uncountable. Write $\hilbertcube_{\lambda}$ for $\lbrack 0,\, 1\rbrack^{\lambda}$.

Assume that ${\mathbf 1}_{{\mathbb B}(\kappa)}\forces ``\dot{G}\subseteq \hilbertcube_{\lambda} \mbox{ is a G$_{\delta}$ set and }\check{X} \subseteq \dot{G}\subseteq \dot{\hilbertcube}_{\lambda}".$ Since the random real forcing notion preserves the Lindel\"of property, we have ${\mathbf 1}_{{\mathbb B}(\kappa)}\forces``\check{X} \mbox{ is Lindel\"of}$". Thus, 
by Lemma \ref{gdeltasets}, we may assume $\dot{G}$ is such that there is a name $\dot{C}$ such that
\[
  {\mathbf 1}_{{\mathbb B}(\kappa)}\forces ``\dot{C} \mbox{ is a countable subset of }\check{\lambda} \mbox{ and } \dot{G} \mbox{ is homeomorphic to } \Pi_{\dot{C}}\lbrack \dot{G}\rbrack \times \lbrack 0,\, 1\rbrack^{\check{\lambda}\setminus \dot{C}}".
\]
Since ${\mathbb B}(\kappa)$ has the countable chain condition, choose a countable set $C\subseteq \lambda$ such that 
${\mathbf 1}_{{\mathbb B}(\kappa)}\forces ``\dot{C} \subseteq\check{C}"$. Then we have 
\[
  {\mathbf 1}_{{\mathbb B}(\kappa)}\forces ``\dot{G} \mbox{ is homeomorphic to } \Pi_{\check{C}}\lbrack \dot{G}\rbrack \times \lbrack 0,\, 1\rbrack^{\check{\lambda}\setminus \check{C}}".
\]
Since in the ground model $\Pi_C\lbrack X\rbrack$ is Hurewicz the first part implies:
\[
  {\mathbf 1}_{{\mathbb B}(\kappa)}\forces ``\Pi_{\check{C}}\lbrack \check{X}\rbrack \subseteq \Pi_{\check{C}}\lbrack \dot{G}\rbrack, \mbox{ a ${\sf G}_{\delta}$ subset of }\lbrack 0,\, 1\rbrack^{\check{C}} \mbox{ and }\Pi_{\check{C}}\lbrack \check{X}\rbrack \mbox{ is Hurewicz}".
\]
Choose a ${\mathbb B}(\kappa)$ name $\dot{F}$ such that 
\[
  {\mathbf 1}_{{\mathbb B}(\kappa)}\forces ``\dot{F}\subseteq \dot{\hilbertcube}^{\check{C}} \mbox{ is an ${\sf F}_{\sigma}$ set and }\Pi_{\check{C}}\lbrack \check{X}\rbrack \subseteq \dot{F} \subseteq \Pi_{\check{C}}\lbrack \dot{G}\rbrack".
\]
But then applying the inverses of the projection maps we find
\[
  {\mathbf 1}_{{\mathbb B}(\kappa)}\forces ``\Pi^{\leftarrow}_{\check{C}}\lbrack\dot{F}\rbrack\subseteq \dot{\hilbertcube}^{\check{\lambda}} \mbox{ is an ${\sf F}_{\sigma}$ set and } \check{X} \subseteq \Pi^{\leftarrow}_{\check{C}}\lbrack\dot{F}\rbrack \subseteq \dot{G}".
\]
It follows that ${\mathbf 1}_{{\mathbb B}(\kappa)}\forces ``\check{X} \mbox{ is a Hurewicz space}."$
$\epf$

\begin{center}{\bf Forcing with countably closed partially ordered sets.}\end{center}

The Tychonoff power $\{0,1\}^{\omega_1}$ of the two-point discrete space is compact and thus has the Hurewicz (and thus Menger) property. But it is not indestructibly Lindel\"of, and so unlike the Rothberger property (Theorem \ref{countablyclosed}), the Hurewicz and Menger properties are not preserved by countably closed forcing. Juh\'asz's spaces (Section 4, Example 4) are Hurewicz spaces but not indestructibly Lindel\"of and so give points ${\sf G}_{\delta}$ examples of this fact. Call a space \emph{indestructibly Hurewicz} (respectively \emph{indestructibly Menger}) if in any generic extension by a countably closed partially ordered set the space is still Hurewicz (respectively Menger). It is clear that indestructibly Hurewicz implies indestructibly Menger, which implies indestructibly Lindel\"of. The converse implications are false. However virtually the same proof as that of Theorem \ref{countablyclosed} gives:
\begin{theorem}\label{indlindhurewicz} The following are equivalent for a space $X$:
\begin{enumerate}
  \item{$X$ is indestructibly Hurewicz.}
  \item{$X$ has the Hurewicz property and is indestructibly Lindel\"of.}  
\end{enumerate}
\end{theorem}

\begin{theorem}\label{indlindmenger} The following are equivalent for a space $X$:
\begin{enumerate}
  \item{$X$ is indestructibly Menger.}
  \item{$X$ has the Menger property and is indestructibly Lindel\"of.}  
\end{enumerate}
\end{theorem}

\begin{center}{\bf When Lindel\"of implies Gerlits-Nagy}\end{center} 

Gerlits and Nagy also introduced the notion of a $\gamma$ space in \cite{GN}: According to \cite{GN} an open cover $\mathcal{U}$ of a space $X$ is an $\omega$-cover if $X\not\in\mathcal{U}$, but for each finite set $F\subseteq X$ there is a $U\in\mathcal{U}$ such that $F\subseteq U$. The symbol $\Omega$ denotes the collection of (open) $\omega$-covers of a space. Call an open cover $\mathcal{U}$ of a space a $\gamma$-cover if it is infinite, and each infinite subset of $\mathcal{U}$ still is a cover of $X$. Let $\Gamma$ denote the collection of (open) $\gamma$-covers of a space.

Then a space is a $\gamma$-space if the selection principle $\sone^{\omega}(\Omega,\Gamma)$ holds: For each sequence $(\mathcal{U}_n:n\in\naturals)$ there is a sequence $(U_n:n\in\naturals)$ such that for each $n$, $U_n\in\mathcal{U}_n$, and each $x\in{\mathbb G}$ is in all but finitely many $U_n$'s. Each $\gamma$-space is a Gerlits-Nagy space. 

A topological space is a \emph{P-space} if each countable intersection of open sets is still an open set. According to \cite{GN} the implication $1 \Rightarrow 4$ in Theorem \ref{PspaceHurewicz} is due to F. Galvin. Since the remainder of these implications are trivial, we attribute this theorem to Galvin: 
\begin{theorem}[Galvin]\label{PspaceHurewicz} Let $X$ be a P-space. The following are equivalent:
\begin{enumerate}
  \item{$X$ is Lindel\"of.}
  \item{$X$ is a Rothberger space.}
  \item{$X$ is a Gerlits-Nagy space.}
  \item{$X$ is a $\gamma$-space.}
\end{enumerate}
\end{theorem}
{\flushleft{\bf Proof:}} It is clear that $4\Rightarrow \, 3\, \Rightarrow\, 2\, \Rightarrow 1$. For the implication $1\Rightarrow 4$, see the Lemma preceding Theorem 3 of \cite{GN}. $\epf$

In particular, Lindel\"of P-spaces are indestructibly Lindel\"of.

L. Zdomskyy \cite{Z} (Theorem 5) proved that it is consistent, relative to the consistency of ZFC, that for regular Lindel\"of spaces the Rothberger property is equivalent to the Gerlits-Nagy property\footnote{Zdomskyy shows that $\mathfrak{u}<\mathfrak{g}$ implies that regular Rothberger spaces have the Hurewicz property. Zdomskyy states that the blanket assumption for the paragraph containing Theorem 5's proof is that the spaces are hereditarily Lindel\"of. This hypothesis is not used in proving Theorem 5, but is used in applying the result of Theorem 5 to an additivity problem. Zdomskyy states his results for paracompact Lindel\"of spaces. It is well-known that regular Lindel\"of spaces are paracompact.}. The Continuum Hypothesis implies the existence of a Lusin set of real numbers. Lusin sets do not have the Gerlits-Nagy property, but have the Rothberger property. 
Among general topological spaces Lusin spaces need not be Rothberger spaces. This can be seen by considering a Sierpi\'nski set of reals in the density topology. But Sierpinski sets have the Menger property in the density topology and, depending on axioms, may have the Hurewicz property in the density topology. For more on these matters see \cite{AT} and \cite{MSierp}.

\begin{center}{\bf Gerlits-Nagy spaces and products.}\end{center}

Next we examine products in the class of Gerlits-Nagy spaces.
Corollary \ref{problem6.6} solves Problem 6.6 of \cite{Ts} and Problems 3.1 through 3.3 of Samet and Tsaban, listed in Section 3 of \cite{SPM18}:
\begin{corollary}\label{problem6.6}
It is consistent that there is a set $X$ of real numbers such that $X$ has the Gerlits-Nagy property while $X^2$ does not, and yet all finite powers of $X$ have the Rothberger property.
\end{corollary}
{\bf Proof:}
Start with a model of CH. Take a Sierpi\'nski set $X$ such that $X\times X$ does not have the Menger property (see for example the remark following Theorem 2.1 of \cite{coc2}). Sierpi\'nski sets are strongly Hurewicz (see e.g. \cite{ST}). Since $X$ is a separable metric space it has the Lindel\"of property in all finite powers.

Now add $\omega_1$ Cohen reals. By Theorem \ref{cohenrothberger}, in the generic extension $X$ has the Rothberger property in all finite powers (thus the Menger property in all finite powers). By Theorem \ref{stronghurewiczpreserve}, in the generic extension $X$ has the Hurewicz property and thus the Gerlits-Nagy property, but by Theorem \ref{nonhurewiczpreserve} $X^2$ does not have the Hurewicz property, and thus not the Gerlits-Nagy property.
$\epf$

\section{Examples }

The following diagram indicates the relationships among the main classes of Lindel\"of spaces considered in this paper:

\begin{figure}[h]
\unitlength=.95mm
\begin{picture}(140.00,60.00)(10,10)

\put(30.00,20.00){\makebox(0,0)[cc]
{\shortstack {Gerlits-Nagy\\ $\{Example\, 2\}$         } }}
\put(70.00,20.00){\makebox(0,0)[cc]
{\shortstack {Rothberger\\ $\{ Lusin\, set   \} $   } }}
\put(110.00,20.00){\makebox(0,0)[cc]
{\shortstack {Ind. Lindel\"of\\ $\{Irrationals   \}$   } }}
\put(30.00,60.00){\makebox(0,0)[cc]
{\shortstack {Hurewicz\\ $\{Example\, 4\} $   } }}
\put(70.00,60.00){\makebox(0,0)[cc]
{\shortstack {Menger\\   . } }}
\put(110.00,60.00){\makebox(0,0)[cc]
{\shortstack {Lindel\"of\\ .  } }}
\put(38.00,61.00){\vector(1,0){20.00}}
\put(80.00,61.00){\vector(1,0){20.00}}
\put(30.00,25.00){\vector(0,1){29.00}}
\put(70.00,25.00){\vector(0,1){29.00}}
\put(110.00,25.00){\vector(0,1){29.00}}
\put(38.00,20.00){\vector(1,0){20.00}}
\put(80.00,20.00){\vector(1,0){20.00}}
\end{picture}
\caption{Classes of points ${\sf G}_{\delta}$ Lindel\"of spaces \label{cshl}}
\end{figure}
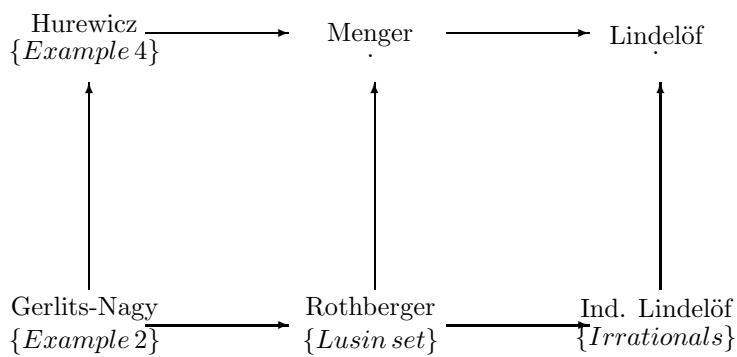

In Figure \ref{cshl} we depict the relationship among six classes of uncountable Lindel\"of spaces with points ${\sf G}_{\delta}$. The spaces in the bottom row are indestructibly Lindel\"of, while the spaces in the top row are not. Example 2 demonstrates that in ZFC the class in the bottom left is nonempty. Example 4 demonstrates that in ZFC the class of Hurewicz spaces are not included in the class of indestructibly Lindel\"of spaces. 

Except for the Rothberger and Gerlits-Nagy classes, the remaining classes are all distinct in ZFC: since the example in the top left corner is not a space as in the bottom right corner, no implication to the top row is an equivalence. Since the space of irrationals is not Menger, no implication from the middle column to the right column is an equivalence. It is consistent that also no implication from the first column is an equivalence: since a Lusin set is not Hurewicz, no implication from the first column to the middle one is an equivalence. For regular Lindel\"of spaces it is independent whether the Rothberger and Gerlits-Nagy classes are distinct. It is true in ZFC that Menger does not imply Hurewicz.

{\flushleft{\bf Example 1}} One-point compactification of an uncountable discrete space:

Let $X$ be an uncountable discrete space and let $X^*$ denote its one-point compactification. Let $\infty$ be the additional point. Then $X^*$ is a Rothberger space but the point $\infty$ is not a ${\sf G}_{\delta}$. Being compact, $X^*$ is also a Hurewicz space, and thus a Gerlits-Nagy space.

{\flushleft{\bf Example 2}} Moore's L-space is hereditarily Rothberger and hereditarily Hurewicz:

In \cite{JTM} Moore solved the famous L-space problem by constructing a non-separable zero-dimensional hereditarily Lindel\"of subspace $X$ of the $\omega_1$-power of the unit circle. This space has among others the following properties:
\begin{enumerate}
  \item{$\vert X\vert = \aleph_1$;}
  \item{$X$ is non-separable;}
  \item{$X^2$ is not Lindel\"of;}
  \item{$X$ remains a hereditarily Lindel\"of space in c.c.c. generic extensions;}
  \item{Every continuous metrizable image of any subspace of $X$ is countable.}      
\end{enumerate}
Item 5 implies that $X$ is hereditarily Rothberger: consider a subspace $Y$ of $X$, which we may assume is infinite. Let $(\mathcal{U}_n:n<\infty)$ be a sequence of open covers of $Y$. Since $Y$ is zero-dimensional and Lindel\"of we may assume each $\mathcal{U}_n$ is countable and disjoint and enumerate it as $(U^n_k:k<\infty)$. Define a function $F:Y\rightarrow\,^{\omega}\omega$ so that for $y\in Y$ we have: $F(y)(m)=k$ if, and only if, $y\in U^m_k$. Then $F$ is continuous and so by item 5, $F[Y]$ is a countable subset of $^{\omega}\omega$. Choose a $g\in\,^{\omega}\omega$ such that for each $y\in Y$ the set $\{n:F(y)(n)=g(n)\}$ is infinite. Then $(U^n_{g(n)}:n<\infty)$ witnesses the Rothberger property of $Y$ for the sequence $(\mathcal{U}_n:n<\infty)$.

Since countable subsets of $^{\omega}\omega$ are bounded, we also see that $X$ has the Hurewicz property hereditarily. Thus, $X$ is hereditarily a Gerlits-Nagy space.   
 
It follows from item 3 that $X$ is a ZFC example of the total failure of preservation of selection properties under finite powers.
Since the space is a Lindel\"of subspace of a power of a Tychonoff space, it is ${\sf T}_4$. And then, since it is hereditarily Lindel\"of, its points are ${\sf G}_{\delta}$: consider an $x\in X$. For each $y\in X\setminus\{x\}$ choose an open neighborhood $U_y$ whose closure does not contain $x$. Then $\{U_y:y\in X\setminus\{x\}\}$ is an open cover of $X$.  Since $X\setminus\{x\}$ is Lindel\"of, choose a countable subset covering $X\setminus\{x\}$. The complements of the closures of the sets in this countable cover witness that $x$ is a ${\sf G}_{\delta}$ point. Moore's spaces are also not first countable (following from item 4) but they are Fr\'echet. These spaces can also be realized as subspaces of $(^{\omega_1}\integers,\mathcal{T})$, where $\mathcal{T}$ is the usual product topology on $\omega_1$ copies of the discrete space of integers. Moreover, as TWO has a winning strategy in $\gone^{\omega_1}(\open,\open)$ on $X$, $X$ remains hereditarily Lindel\"of in every countably closed forcing extension. 

{\flushleft{\bf Example 3}} Gorelic's example:

Gorelic's space, ${\mathbb G}$, is obtained by countably closed forcing ${\mathbb P}$. An uncountable cardinal number $\kappa$ with $\kappa^{\aleph_1}=\kappa$ is taken, and then a subspace ${\mathbb G}\subseteq \,^{\kappa}\{0,1\}$ of cardinality $\kappa$ is obtained in a generic extension which has $2^{\aleph_1}=\kappa$, such that ${\mathbb G}$ is zero-dimensional, Lindel\"of, and has each point ${\sf G}_{\delta}$.

Indeed, ${\mathbb G}$ has both the Rothberger property and the Hurewicz property, as we shall now show. Towards this we first describe Gorelic's construction: first, choose a partition of $\kappa$ into pairwise disjoint countable sets $A_{\alpha}$, $\alpha<\kappa$. And for each $\beta<\kappa$ choose a countable dense subset $H_{\beta}$ of $^{A_{\beta}}\{0,1\}$ and enumerate it as $\{h^{\beta}_{\alpha}:\alpha\in A_{\beta}\}$.

The elements of ${\mathbb P}$ are pairs $(F,T)$ where $F$ and $T$ are as follows: 
\begin{itemize}
  \item {Both $F$ and $T$ are countable sets;}
  \item {There is a countable set $C\subseteq\kappa$ such that each element of $F$ is a function with domain $D = \bigcup_{\beta\in C}A_{\beta}$;}
  \item {$F$ can be indexed as $\{f_{\xi}:\xi\in D\}$ such that whenever $\beta\in C$ and $\alpha\in A_{\beta}$ then $f_{\alpha}\lceil A_{\beta} = h^{\beta}_{\alpha}$;}
  \item {Each element of $T$ is a set $B\subseteq {\sf Fn}(D,2)$, with $\{[\sigma]\subseteq \,^{A}\{0,1\}:\sigma\in B\}$ a cover of $F$.}
\end{itemize} 

The order on ${\mathbb P}$ is defined as follows:
\[
  q = (F_2,T_2) \le (F_1,T_1) = p
\]
if, letting $C_p$ and $C_q$ be the associated $C$'s, and $D_p$ and $D_q$ the associated $D$'s, we have:
\begin{enumerate}
 \item{$C_q\supseteq C_p$;}
 \item{For each $\alpha\in D_p$, $f^q_{\alpha}\lceil_{D_p} = f^p_{\alpha}$;}
 \item{$T_q\supseteq T_p$.} 
\end{enumerate}

When $H$ is a ${\mathbb P}$-generic filter, put for each $\alpha\in\kappa$, 
\[
   f_{\alpha} = \bigcup\{f^p_{\alpha}:p\in H \mbox{ and } \alpha\in D_p\}.
\]
 Then put ${\mathbb G}=\{f_{\alpha}:\alpha\in\kappa\}$. Let $\Gamma$ be a ${\mathbb P}$-name for ${\mathbb G}$.

Here are some observations about $p=(F,T)\in{\mathbb P}$ (in the above notation): 

\begin{enumerate}
  \item {If $U\in T$ then $p\forces ``(\forall f\in\Gamma)(\exists \sigma\in U)(\sigma\subseteq f)"$.}
  \item {If $\tau$ is a ${\mathbb P}$-name such that $p\forces ``\tau \mbox{ is a cover of $\Gamma$ by basic open sets}"$ then there is a $q<p$ with a $B\in T_q$ such that $q\forces ``\{[\sigma]\subseteq \,^{\kappa}\{0,1\}: \sigma\in B\} \subseteq \tau \mbox{ covers }\Gamma"$.}
  \item {For $B\subseteq {\sf Fn}(D_p,2)$ such that $\{[\sigma]\subseteq \,^{D_p}\{0,1\}: \sigma\in B\}$ covers $F\subseteq \,^{D_p}\{0,1\}$ put $q = (F,T\cup\{B\})$. Then $q\in{\mathbb P}$ and $q\le p$.}
\end{enumerate}
By the first item above, $p$ forces that $\{[\sigma]\subseteq \,^{\kappa}\{0,1\}: \sigma\in B\}$ is an open cover of ${\mathbb G}$ whenever $B\in T_p$. 

To see that for each $p\in{\mathbb P}$ we have
\[
  p \forces ``\Gamma \mbox{ has the Rothberger property}"
\] 
consider ${\mathbb P}$-names $\tau_n$ such that 
\[
  p \forces ``(\forall n)(\tau_n \mbox{ is a basic open cover of }\Gamma)"
\]
Now using the second remark above, recursively choose a sequence $p_n = (F_n,T_n)$ and $B_n$, $n<\infty$, such that:
\begin{enumerate}
  \item {$B_n\in T_n$ and $p_n\forces ``\{[\sigma]\subseteq \,^{\kappa}\{0,1\}: \sigma\in B_n\} \subseteq \tau_n \mbox{ covers }\Gamma"$.}
  \item {$p_1< p$ and for all $n$, $p_{n+1}<p_n$.}
\end{enumerate}
Put $D = \bigcup_{n<\infty}D_{p_n}$, a countable subset of $\kappa$, and for $\alpha\in D$ put $f_{\alpha} = \cup\{f^{p_n}_{\alpha}:\alpha\in D_{p_n}\}$. Put $F=\{f_{\alpha}:\alpha\in D\}$ and $T=\bigcup\{T_n:n<\infty\}$. Then $q=(F,T)$ is an element of ${\mathbb P}$ and for all $n$, $q\le p_n$.

Now $F\subseteq \,^{D}\{0,1\}$ is countable, and for each $n$, $\mathcal{U}_n = \{[\sigma]\subseteq\,^D\{0,1\}:\sigma\in B_n\}$ is an open cover of $F$. Since $F$ is countable we can choose for each $n$ a $\sigma_n\in B_n$ such that $\{[\sigma_n]\subseteq\,^{D}\{0,1\}:n<\infty\}$ is an open cover of $F$. Note that $B = \{\sigma_n:n<\infty\}\subseteq {\sf Fn}(D,2)$. But then we have
\[
  p_{\omega} = (F,T\cup\{B\}) \in{\mathbb P}.
\]
Then $B$ witnesses\footnote{Define the ${\mathbb P}$-name $\mu=\{(\langle n,\sigma_n\rangle,{\mathbf 1}):n<\infty\}$. Then check that $p_{\omega}$ forces that $\mu$ codes a sequence of such $U_n$'s.} that
\[
  p_{\omega}\forces ``(\forall n)(\exists U_n\in\tau_n)(\{U_n:n<\infty\} \mbox{ covers }\Gamma)".
\]
We have shown that the set of conditions forcing that ${\mathbb G}$ has the Rothberger property is dense below any element of ${\mathbb P}$.

To show that ${\mathbb G}$ has the Hurewicz property, instead of choosing singletons from each $B_n$ choose finite sets such that 
each element of $F$ is covered by all but finitely many of the unions of these finite sets, and let the $B$ in the definition of $p_{\omega}$ be the union of these finite sets.

With a little bit more work one can show that ${\mathbb G}$ is a $\gamma$-space. It follows that ${\mathbb G}$ is a $\gamma$-space with points ${\sf G}_{\delta}$ and yet has cardinality $2^{\aleph_1}>\aleph_1 = 2^{\aleph_0}$.

{\flushleft{\bf Example 4}} Juh\'asz's examples:

In \cite{Juhasz}, Example 7.2, Juh\'asz gave examples of Lindel\"of spaces with points ${\sf G}_{\delta}$ which are of arbitrary large cardinality below the first measurable cardinal. We give, for the reader's convenience, the description from \cite{FDT} for these. 

Consider an infinite cardinal $\kappa$ less than the first measurable cardinal. Recursively define\\
\begin{tabular}{lll}
$\kappa_0$ & $=$ & $\kappa$\\
$\kappa_{n+1}$ & $=$ & $2^{\kappa_n}$\\
$\kappa_{\omega}$ & $=$ & $\sup\{\kappa_n:n<\infty\}$\\
\end{tabular}

For an infinite set $S$, let $\mathcal{F}(S)$ denote the set of non-principal ultrafilters on $S$. For $A\subseteq S$ put \^{A} $= \{{\mathcal U}\in\mathcal{F}(S):A\in\mathcal{U}\}$. And for $\mathcal{U}\in{\mathcal F}({\mathcal F}(S))$ put $\mathcal{U}^{\prime} = \{A\subseteq S: \mbox{\^{A}}\in\mathcal{U}\}$. Then $\mathcal{U}^{\prime}$ is an ultrafilter and a member of $\mathcal{F}(S)$. 

Define $J_0=\kappa$ and for each $n$ define $J_{n+1} = \mathcal{F}(J_n)$ Finally, put ${\mathbb J}_{\kappa} = \bigcup_{n<\infty}J_n$.

Now for each $n$, and for $\mathcal{U}\in X_{n+1}$ define $\mathcal{U}^{(i)}$ for $i\le n$ as follows:
\[
  \mathcal{U}^{(0)} =\mathcal{U},\,\, \mathcal{U}^{(i+1)} = (\mathcal{U}^{(i)})^{\prime}.
\]
The topology on ${\mathbb J}_{\kappa}$ is defined by describing neighborhood bases for points:
\begin{itemize}
  \item {Points of $J_0$ are isolated;}
  \item {For $\mathcal{U}\in J_{n+1}$, all sets of the form
   \[
     \{{\mathcal U}\}\cup(\bigcup_{i=0}^n A^{(i)}), \,\, A^{(i)}\in\mathcal{U}^{(i)},\, 0\le i\le n.
   \]
        constitute a neighborhood base of $\mathcal{U}$.}
\end{itemize}
These spaces are not ${\sf T}_2$. Juh\'asz showed in \cite{Juhasz} that these spaces have each point ${\sf G}_{\delta}$. Tall showed in \cite{FDT} Theorem 19 that these spaces are destructibly Lindel\"of. Thus, they are not Rothberger.  

The following crucial observation of Juh\'asz provides a tool to show that spaces of the form ${\mathbb J}_{\kappa}$ are Hurewicz:

\begin{lemma}[Juh\'asz]\label{juhaszlemma} For any set $S$ and function $f:{\mathcal F}(S)\longrightarrow \mathcal{P}(S)$ such that for each $\mathcal{U}\in\mathcal{F}(S)$ we have $f(\mathcal{U})\in\mathcal{U}$, there are finitely many ultrafilters $\mathcal{U}_1,\cdots\mathcal{U}_n$ in $\mathcal{F}(S)$ such that $\vert S\setminus \bigcup_{1\le i\le n} f(\mathcal{U}_i)\vert<\aleph_0$. 
\end{lemma}
Lemma \ref{juhaszlemma} implies for each $n$ that from each open cover of $J_{n+1}$ one can choose finitely many members whose union covers all but finitely many points of $J_n$.

We can show that the property of this Lemma is preserved in the random real extension, and thus ${\mathbb J}_{\kappa}$ remains a Hurewicz space in the random real extension. We don't know if this is also the case with Cohen reals.

\begin{lemma}\label{Juhaszspacehrewicz}
For each infinite cardinal $\kappa$ the space ${\mathbb J}_{\kappa}$ has the Hurewicz property.
\end{lemma}
{\bf Proof:} Apply Lemma \ref{juhaszlemma}:
given a sequence $({\mathcal G}_n:n<\infty)$ of open covers of ${\mathbb J}_{\kappa}$, choose for each $n$ a finite set $\mathcal{H}_n\subseteq \mathcal{G}_n$ which covers $\cup_{i\le n}J_i$. It follows that each element of $J_{\omega}$ is contained in all but finitely many $\cup\mathcal{H}_n$.
$\epf$

Another interesting property of these spaces is given in Theorem 21 of \cite{BT}:
\begin{theorem}[Baumgartner-Tall]\label{subspaces1} ${\mathbb J}_{\kappa}$ does not have a Lindel\"of subspace of cardinality $\aleph_1$.
\end{theorem}

Consider ${\mathbb J}_{\kappa}$ from the ground model in a generic extension by ${\mathbb P}(\aleph_1)$, the partially ordered set for adding $\aleph_1$ Cohen reals. By 
 Theorem \ref{cohenrothberger} ${\mathbb J}_{\kappa}$ acquired the Rothberger property.
  Now recall Lemma 25 of \cite{BT}:
\begin{lemma}[Baumgartner-Tall]\label{nonlindelofpreserve} Suppose X is a space in the ground model.  Force with a Property K partial order. Then if X has a Lindel\"of subspace of size $\aleph_1$ in the extension, it has one in the ground model.
\end{lemma}
Both the Cohen reals partial order ${\mathbb P}(\kappa)$ and the random reals partial order ${\mathbb B}(\kappa)$ have property K. 

It follows that ${\mathbb J}_{\kappa}$ still does not have a Lindel\"of subspace of cardinality $\aleph_1$ in the generic extension. Thus, it is consistent relative to the consistency of ZFC that there are spaces of arbitrary large cardinality below the first measurable cardinal which have points ${\sf G}_{\delta}$ and the Rothberger  
 property, and no Lindel\"of subspace of cardinality $\aleph_1$

Examining Corollary 24 of \cite{BT} we see that it is consistent that there is a ${\sf T}_1$ Rothberger 
space of cardinality $\aleph_2=(2^{\aleph_0})^+$ which has no Lindel\"of subspace of cardinality $\aleph_1$. Indeed: if we take a ${\mathbb J}_{\kappa}$, add $\aleph_1$ Cohen reals to make it a Rothberger  
 space, and then collapse its cardinality to $\aleph_2$ with conditions of size $\le \aleph_1$, we obtain a Rothberger 
  space with no Lindelof subspace of size $\aleph_1$. 

{\flushleft{\bf Example 5}} Koszmider and Tall's example:

Juh\'asz's spaces are not ${\sf T}_2$ but have points ${\sf G}_{\delta}$. In \cite{KT} Koszmider and Tall constructed by countably closed forcing ${\mathbb P}$ over a model of CH a subspace ${\mathbb K}$ of $(^{\omega_2}2,\mathcal{T}_2)$. Thus ${\mathbb K}$ is a P-space which is ${\sf T}_3$ and Lindel\"of (thus ${\sf T}_4$), is of cardinality $\aleph_2$, and it has no Lindel\"of subspace of cardinality $\aleph_1$. It follows from Theorem \ref{PspaceHurewicz} that ${\mathbb K}$ has the Gerlits-Nagy property and yet fails to have a Lindel\"of subspace of cardinality $\aleph_1$. 

{\flushleft{\bf Example 6}} Stationary Aronszajn lines:

In \cite{stevo} Todor\v{c}evi\'c investigates the status of certain selection principles on the class of Aronszajn lines. These are first countable spaces, and in Theorem 1 of \cite{stevo} it is shown that stationary Aronszajn lines are Rothberger spaces. Theorem 3 of \cite{stevo}, combined with a result 
 of \cite{BCM} that for each Lindel\"of space $X$, if each continuous function from $X$ to $^{\omega}\reals$ is bounded in the eventual domination order, then the space $X$ is a Hurewicz space, shows these also have the Hurewicz property. Thus, the stationary Aronszajn lines have the Gerlits-Nagy property. In Theorem 6 of \cite{stevo} Todor\v{c}evi\'c shows that the more restrictive class of stationary Countryman lines indeed are $\gamma$-spaces, and that though this property is preserved by finite powers, it is not preserved by finite unions, and not even the Lindel\"of property is preserved by finite products of $\gamma$-spaces. Since the Rothberger property and the Hurewicz property are preserved by countable unions, these examples also demonstrate that the class of $\gamma$ spaces is more restrictive than the class of Gerlits-Nagy spaces. Each of the mentioned examples can be taken to be a space of cardinality $\aleph_1$.

\end{document}